\documentclass[11pt,a4paper,reqno]{amsart}

\usepackage{graphicx,mathrsfs,amssymb,amsmath,mathtools,bm}
\usepackage[dvipsnames]{xcolor}
\usepackage[colorlinks=true,linkcolor=blue,citecolor=red,urlcolor=magenta]{hyperref}
\usepackage{enumerate}

\usepackage[backgroundcolor=yellow, colorinlistoftodos,prependcaption,textsize=small,textwidth=25mm]{todonotes}
\usepackage{cancel}
\usepackage[normalem]{ulem}

\newtheorem{theorem}{Theorem}[section]
\newtheorem{lemma}[theorem]{Lemma}
\newtheorem{corollary}[theorem]{Corollary}
\newtheorem{proposition}[theorem]{Proposition}
\theoremstyle{definition}
\newtheorem{definition}[theorem]{Definition}
\newtheorem{remark}[theorem]{Remark}

\numberwithin{equation}{section}

\newcounter{assumption}
\renewcommand{\theassumption}{\Alph{assumption}}

\newcommand{\C}{{\mathbb C}}
\newcommand{\Z}{{\mathbb Z}}
\newcommand{\T}{{\mathbb T}}

\newcommand{\N}{{\mathbb N}}
\newcommand{\R}{{\mathbb R}}
\newcommand{\bk}{\mathbf{k}}
\newcommand{\bde}{\bm{\delta}}
\newcommand{\bh}{\mathbf{h}}

\newcommand{\FT}{\mathscr{F}_{_{\T^n}}}

\newcommand{\CTE}{C^\infty(\T^n,E)}

\newcommand{\WTE}{W_p^k(\T^n,E)}

\newcommand{\dbar}{{\textstyle d^{\hspace*{-0.4em}-}\hspace*{-0.15em}}}

\newcommand{\sumint}{\mathop{\mbox{$\displaystyle\sum\limits_{\bk\in\Z^n}\int_{\T^n}$}}}
\newcommand{\ossumint}{\mathop{\mbox{$\displaystyle\mathrm{Os}\!-\!\!\!\!\sum\limits_{\bk\in\Z^n}\int_{\T^n}$}}}

\newcommand{\op}{\mathop{\textrm{\upshape op}}}

\renewcommand{\Re}{\mathop{\mathrm{Re}}}
\renewcommand{\Im}{\mathop{\mathrm{Im}}}

\renewcommand{\epsilon}{\varepsilon}
\renewcommand{\bar}[1]{\overline{#1}}
\renewcommand{\tilde}[1]{\widetilde{#1}}

\allowdisplaybreaks
\begin{document}
\color{black}

\title[Pseudodifferential operators on the torus]{Generation of semigroups for vector-valued pseudodifferential operators on the torus with space-dependent symbols}
\author{B. Barraza Mart\'{i}nez}
\address{B. Barraza Mart\'{i}nez, Universidad del Norte, Departamento de Matem\'aticas, F\'isica y Ciencia de Datos,  Barranquilla (Colombia)}
\email{bbarraza@uninorte.edu.co}

\author{R.\ Denk}
\address{R.\ Denk, Universit\"at Konstanz, Fachbereich f\"ur Mathematik und Statistik, Konstanz (Germany)}
\email{robert.denk@uni-konstanz.de}

\author{J. Hern\'andez Monz\'on}
\address{J. Hern\'andez Monz\'on, Universidad del Norte, Departamento de Matem\'aticas, F\'isica y Ciencia de Datos,  Barranquilla (Colombia)}
\email{jahernan@uninorte.edu.co}

\author{M.\ Nendel}
\address{M.\ Nendel, University of Waterloo, Department of Statistics and Actuarial Science, Waterloo, Ontario (Canada)}
\email{mnendel@uwaterloo.ca}

\date{\today}

\thanks{The authors would like to thank COLCIENCIAS (Project 121556933488) and DAAD
for financial support.\ The last-named author gratefully acknowledges financial support through the Deutsche Forschungsgemeinschaft (DFG, German Research Foundation) -- SFB 1283/2 2021 -- 317210226 and the Natural Sciences and Engineering Research Council of Canada via Discovery Grant no.\ RGPIN-2025-04219.}

\begin{abstract}
We consider toroidal pseudodifferential operators with spatially
dependent, operator-valued symbols of positive order and finite regularity.
For an arbitrary complex Banach space $E$, we prove that the negative of
the closure of a parabolic operator on $C^\infty(\T^n,E)$ generates an
analytic semigroup on $W_p^k(\T^n,E)$,
where $1\le p<\infty$ and $k\in\N_0$, provided the symbol has sufficiently
many spatial derivatives and discrete differences. The proof combines
uniform estimates for parameter-dependent kernels with right and left
parametrices. It does not require the UMD property of $E$.
Real interpolation gives corresponding resolvent estimates on
$B_{pq}^s(\T^n,E)$ for $0<s<k$ and $1\le q\le\infty$.
The induced semigroups are strongly continuous for $q<\infty$ and on
the little Besov spaces when $q=\infty$.
\end{abstract}

\let\BDHNsubjclassname\subjclassname
\renewcommand{\subjclassname}{{\BDHNsubjclassname}}

\maketitle

\section{Introduction}

 In this article, we study vector-valued  pseudodifferential operators on the torus $\T^n = (\R/2\pi \Z)^n$ and  generation of semigroups in Sobolev and Besov spaces. The Fourier series approach to toroidal pseudodifferential
operators allows a global quantization with covariable $\bk\in\Z^n$, and we define such operators
via their symbols, which are functions of the space variable $x\in\T^n$ and the covariable $\bk\in\Z^n$. As a standard reference for toroidal quantization, we mention the monograph of Ruzhansky and Turunen \cite{ruzhansky-turunen10}. The theory of toroidal pseudodifferential operators was developed by, e.g.,  Agranovich \cite{agranovich79}, Amosov \cite{amosov88}, McLean \cite{mclean91}, Melo \cite{melo97}, Saranen--Wendland \cite{saranen-wendland87},  Ruzhansky--Turunen \cite{ruzhansky-turunen07, ruzhansky-turunen09}, Turunen \cite{turunen00}, and Turunen--Vainikko \cite{turunen-vainikko98}. For mapping properties of toroidal pseudodifferential operators in $L^p$-spaces, we refer to Delgado \cite{delgado13}, Molahajloo--Wong \cite{molahajloo-shahla-wong10}, Wong \cite{wong11}, Cardona \cite{cardona14} and others.
Mapping properties of scalar pseudodifferential 
operators in non-reflexive end-point spaces were also studied in Benhamoud \cite{Benhamoud2024Endpoint}. 
While 
the above references deal with scalar symbols,  we study operator-valued symbols and 
Banach space valued functions. 

Vector-valued functions appear in applications, for instance, if the unknown function depends on an additional parameter and therefore takes values in a function space in that parameter, see, e.g., coagulation-fragmentation problems in Amann--Walker \cite{amann-walker05}.
Other typical applications concern equations in cylindrical domains with periodic boundary 
conditions as considered, e.g., in Denk--Nau \cite{denk-nau13}, Favini--Guidetti-Yakubov \cite{favini-giudetti-yakubov11}, Nau--Saal \cite{nau-saal12}, and Rabinovich \cite{rabinovich11}.
Results on the nuclearity of operator-valued pseudodifferential operators in Hilbert spaces can be found in Cardona--Kumar \cite{CardonaKumar2021Nuclear}. Further results on mapping properties can be found in Xia--Xiong \cite{XiaXiong2019Mapping}. For applications to partial differential equations on the torus (in the scalar situation), we refer to, e.g., Cardona--Delgado--Ruzhansky \cite{CardonaDelgadoRuzhansky2022Drift, CardonaDelgadoRuzhansky2024Hyperbolic} 
and Guinazu--Vergara \cite{Guinazu-Vergara21}.

In particular, if the unknown function is a density with respect to an additional parameter,
a natural state space is an $L^1$-space with respect to this parameter. This excludes 
the application of results by Arendt--Bu \cite{arendt-bu02, arendt-bu04} and Bu--Kim \cite{bu-kim04}  
on vector-valued Fourier multipliers in $L^p$-spaces, where the state space is assumed to be a UMD space. 
 For more results on Fourier multipliers in UMD spaces
and the transference principle from Fourier transform to Fourier series in such spaces, 
see also Hyt\"onen et al.\ \cite[Chapter 5]{HNVW16} and 
Fackler--Hyt\"onen--Lindemulder \cite{FacklerHytonenLindemulder2020Weighted}. As we want to include $L^1$-spaces, these results cannot be 
applied in the present situation. 

For toroidal pseudodifferential operators  in an arbitrary Banach space, the generation
of holomorphic semigroups was studied in the previous work 
\cite{barraza-denk-hernandez-nau}. There, we considered symbols which were independent of
the space variable $x$. This simplifies the analysis, in particular for the composition of
operators or for the adjoint operator. In pseudodifferential theory, the symbols are mostly
assumed to be infinitely smooth in the space variable. However, from an application point of
view, this is too restrictive. Therefore, in the present paper, we assume only finite smoothness
with respect to $x$, which leads to technical difficulties. We try to give reasonable sufficient
regularity  for each of the results, without claiming them to be sharp. We also want to mention
that a different possible approach uses regularization, cf. Taylor \cite{Taylor91} and Abels \cite{Abels-Pfeuffer18} for the scalar case. However, for this one has to study more general symbol classes than the standard
H\"ormander class $S^{\mu}_{1,0}$ we consider. To summarize, we address two main challenges
in our analysis of toroidal pseudodifferential operators:
\begin{itemize}
  \item The symbols $a(x,\bk)$ have values in the space of bounded linear operators on an arbitrary Banach space $E$, without assuming 
  the UMD property for $E$,
  \item the symbols depend on the space variable $x$, but they are only of finite smoothness with respect to $x$.
\end{itemize}
In this general situation, Theorem~\ref{Theorem2.13} shows that, under appropriate ellipticity assumptions, $-A_k$
generates an analytic $C_0$-semigroup on $W_p^k(\T^n,E)$, where  $A_k$ is the $W_p^k(\T^n,E)$-closure of the pseudodifferential operator defined on $C^\infty(\T^n,E)$. Using real interpolation, this implies related results on Besov spaces, see Corollary~\ref{cor:Besov-generation}. 
For mapping properties on toroidal Besov spaces, we also refer to the previous work \cite{BDHN2018}. 
As we cannot apply Fourier multiplier techniques for the reasons explained above, all proofs
are based on explicit kernel estimates. 

In Section 2 of this paper, we define the class of considered pseudodifferential operators,
using an oscillatory sum-integral, and give results on the composition of such operators,
taking the finite space regularity into account. Section 3 contains the main results on 
the generation of analytic semigroups in Sobolev and Besov spaces, Theorem~\ref{Theorem2.13} and
Corollary~\ref{cor:Besov-generation}. The proof is based on kernel estimates, see, in particular,
Lemma~\ref{lem:remainder-estimate}, which can be shown by an explicit description of the kernel
(Lemma~\ref{Lemma2.9}) and estimates on infinite series, which are uniform with respect to 
the resolvent parameter $\lambda$ (Lemma~\ref{Lemma2.12}). On a technical level, we have to 
deal with parameter-dependent and operator-valued double symbols $p_\lambda(x,y,\bf k)$, see
the proof of Theorem~\ref{Theorem2.13}. By a careful analysis of these symbols, we can show that
the resolvent $(\lambda+A_k)^{-1}$ exists and decays (after a shift) as $1/|\lambda|$, which 
shows that the operator $-A_k$ is sectorial and therefore generates an analytic semigroup.

\section{Preliminaries}
Throughout this paper, we fix $n\in\N$ and a complex Banach space $E$ with norm $\Vert \cdot\Vert$. We consider pseudodifferential operators on the $n$-dimensional torus $\T^n:=(\R/2\pi \Z)^n$ corresponding to symbols with values in the space $L(E)$ of bounded linear operators on $E$. We use $[-\pi,\pi)^n$ as a set of representatives of $\T^n$, and refer to \cite{barraza-denk-hernandez-nau} for the definitions of spaces of functions and distributions on $\T^n$ as well as on the dual space $\Z^n$ of $\T^n$, together with their topologies. For $f\in L^1(\T^n,E)$ the \emph{toroidal Fourier transform} is given by
\begin{equation*}
    (\FT f)(\bk) := \hat{f}(\bk) := \int_{\T^n} e^{-i\bk\cdot x}f(x)\,\dbar x \quad (\bk\in\Z^n),
\end{equation*}
where $\dbar x := (2\pi)^{-n} dx$. The Fourier transform is extended by duality to the space of vector-valued toroidal distributions $\mathscr{D}'(\T^n,E):=L(C^\infty(\T^n), E)$, see \cite[Section 2]{barraza-denk-hernandez-nau} for more details.\\
Now, let $j\in\{1,\dots,n\}$, and let $\bde_j:=(\delta_{jk})_{k=1,\dots,n}$ be the $j$-th unit vector in $\R^n$, where $\delta_{jk}$ is the standard \emph{Kronecker delta} given by $\delta_{jk}=1$ if $k=j$, and $\delta_{jk}=0$ if $k\neq j$. For $a:\Z^n\to E$ and $\alpha\in\N_0^n$, the operators for \emph{discrete derivatives} or \emph{differences} $\Delta_\bk^\alpha$ are defined in the following way:
\begin{align*}
    \Delta_{k_j}a(\bk) & := a(\bk + \bde_j) - a(\bk) \quad (\bk\in\Z^n),\\
    \Delta_\bk^\alpha & := \Delta_{k_1}^{\alpha_1} \cdots \Delta_{k_n}^{\alpha_n}.
\end{align*}
Furthermore, we have the discrete Leibniz formula
\begin{equation}
    \label{4-0a}
\Delta_\bk^\alpha (fg)(\bk) = \sum_{\substack{\beta\in{\N_0^n}\\ \beta\le\alpha}} \binom \alpha\beta \big[ \Delta_\bk^\beta f\big](\bk) \big[\Delta_\bk^{\alpha-\beta} g\big](\bk+\beta)\quad (\bk\in\Z^n)
 \end{equation}
for $f,g\colon\Z^n\to L(E)$, which can be found (in the scalar case) in \cite[Proposition~3.1]{ruzhansky-turunen10a}, see also \cite[Section 3.3.1]{ruzhansky-turunen10} for more details about discrete analysis on the torus. As usual, we set $\langle x\rangle:=(1+\vert x\vert^2)^{1/2}$ for $x\in\R^n$, where $\vert x\vert$ is the Euclidean norm of $x$ in $\R^n$.

\begin{definition}\label{Def_Symbols_Operators}
Let $m\in\R$, $\rho\in\N_0$, and $r\in\N_0$.
    \begin{itemize}
        \item[a)] The symbol class $S^{m,\rho,r}:=S^{m,\rho,r}(\T^n\times \Z^n, L(E))$ consists of all functions $a:\T^n\times \Z^n \to L(E)$ such that $[x\mapsto a(x,\bk)]\in C^r(\T^n,L(E))$ for all $\bk\in\Z^n$, and
        \begin{equation*}
            \Vert a\Vert_m^{(\rho,r)}:= \max_{\vert \alpha\vert\leq \rho}\,\max_{\vert \beta\vert\leq r}\,\sup_{x\in\T^n}\,\sup_{\bk\in\Z^n} \langle \bk\rangle^{\vert \alpha\vert - m}\Vert \Delta_\bk^\alpha\partial_x^\beta a(x,\bk)\Vert_{L(E)}<\infty.
        \end{equation*}
        \item[b)] For $a\in S^{m,\rho,r}$ the pseudodifferential operator $\op[a]$ is defined by
        \begin{equation}\label{Eq_def_op_a}
            (\op[a]f)(x) := \sum\limits_{\bk\in\Z^n}e^{i\bk\cdot x}a(x,\bk)\hat{f}(\bk) \quad \big(f\in C^\infty(\T^n,E), \ x\in\T^n\big).
        \end{equation}
    \end{itemize}
\end{definition}

\begin{remark}\label{Remark_osc-sum-int}
Formally, for $a\in S^{m,\rho,r}$, $f\in C^\infty(\T^n,E)$, and $x\in\T^n$, we have
\begin{align*}
  (\op[a]f)(x) & = \sum\limits_{\bk\in\Z^n}e^{i\bk\cdot x}a(x,\bk)\hat{f}(\bk)\\
  & = \sum\limits_{\bk\in\Z^n}e^{i\bk\cdot x}a(x,\bk)\int_{\T^n} e^{-i\bk\cdot y}f(y)\,\dbar y\\
  & = \sum\limits_{\bk\in\Z^n}\int_{\T^n} e^{i\bk\cdot(x- y)}a(x,\bk)f(y)\,\dbar y.
\end{align*}
While the sum of these integrals converges absolutely, the sum of the integrals of the norms of the integrands need not be finite. We therefore also interpret this expression as an oscillatory sum-integral. More precisely, we write
\begin{align*}
    (\op[a]f)(x) & \coloneqq \ossumint e^{i\bk\cdot(x- y)}a(x,\bk)f(y)\,\dbar y\\
    & \coloneqq \lim\limits_{\varepsilon\searrow 0}\sumint e^{i\bk\cdot(x- y)}\chi(\varepsilon\bk)a(x,\bk)f(y)\,\dbar y,
\end{align*}
where $\chi\in\mathscr{S}(\R^n)$ satisfies $0\leq \chi\leq 1$ and $\chi(0)=1$. This limit does not depend on the choice of $\chi$. Indeed, let
\begin{equation}\label{Eq-L_ky}
 L_{\bk,y} :=\langle \bk\rangle^{-2}(1+i\bk\cdot\nabla_y)\qquad (\bk\in\Z^n,\  y\in\T^n).
\end{equation}
Note that $L_{\bk,y}^Ne^{i\bk\cdot(x-y)}=e^{i\bk\cdot(x-y)}$ for all $N\in\N$. Integration by parts yields
\begin{align*}
  \ossumint & e^{i\bk\cdot(x- y)}a(x,\bk)f(y)\,\dbar y\\
  & = \lim\limits_{\varepsilon\searrow 0}\sumint \big( L_{\bk,y}^N e^{i\bk\cdot(x- y)}\big)\chi(\varepsilon\bk)a(x,\bk)f(y)\,\dbar y\\
  & = \lim\limits_{\varepsilon\searrow 0}\sumint  e^{i\bk\cdot(x- y)}\chi(\varepsilon\bk)a(x,\bk)(L_{\bk,y}^t)^Nf(y)\,\dbar y,
\end{align*}
where
\begin{equation}\label{Eq-L_kyt}
    L_{\bk,y}^t\coloneqq \langle \bk\rangle^{-2}(1-i\bk\cdot\nabla_y) \qquad (\bk\in\Z^n,\  y\in\T^n).
\end{equation}
A straightforward calculation shows that
\begin{equation}\label{Eq-L_kytN}
   (L_{\bk,y}^t)^N = \sum\limits_{\substack{\gamma\in \N_0^n\\   \vert\gamma\vert\leq N}} c_{\gamma,N}(\bk)\partial_y^\gamma  \quad \text{with}\quad \vert c_{\gamma,N}(\bk)\vert \leq C_N\langle \bk\rangle^{-N}\quad (\vert \gamma\vert\leq N).
\end{equation}
Therefore, choosing $N>m+n$, we obtain
\begin{equation*}
  \big\Vert e^{i\bk\cdot(x-y)}\chi(\varepsilon\bk)a(x,\bk)(L_{\bk,y}^t)^Nf(y)\big\Vert \le C_N\Vert a\Vert_m^{(\rho,0)}
  \Vert f\Vert_{C^N(\T^n,E)}\langle\bk\rangle^{m-N}.
\end{equation*}
The right-hand side is summable over $\bk\in\Z^n$ and independent of
$\varepsilon$ and $x,y$. Dominated convergence therefore gives
\[
 (\op[a]f)(x)=\sumint e^{i\bk\cdot(x-y)}a(x,\bk)
                  (L_{\bk,y}^t)^Nf(y)\,\dbar y.
\]
This proves that
the limit exists and is independent of $\chi$.

\end{remark}

For double symbols, $S^{m,\rho,r}(\T^{2n}\times\Z^n,L(E))$ is defined
in the same way, i.e., $a\in S^{m,\rho,r}(\T^{2n}\times\Z^n,L(E))$ if  $a(\cdot,\cdot,\bk)\in C^r(\T^{2n},L(E))$ for every
$\bk\in\Z^n$, and
\[
 \Vert a\Vert_m^{(\rho,r)}:=
 \max_{|\alpha|\le\rho}\max_{|\beta|+|\gamma|\le r}
 \sup_{x,y\in\T^n}\sup_{\bk\in\Z^n}
 \langle\bk\rangle^{|\alpha|-m}
 \Vert\Delta_\bk^\alpha\partial_x^\beta\partial_y^\gamma
              a(x,y,\bk)\Vert_{L(E)}<\infty.
\]

\begin{definition}\label{Def_op_a_xyk}
 Let $m\in\R$, $\rho\in \N$, $r\in\N_0$ with $r\geq \lfloor m  \rfloor + n + 1$, and $\chi\in\mathscr{S}(\R^n)$ with $0\leq\chi\leq1$ and $\chi(0)=1$. For $a\in S^{m,\rho,r}\big( \T^{2n}\times\Z^n , L(E)\big)$  and $f\in C^\infty(\T^n,E)$ we define
 \begin{align}
 \begin{split}\label{Eq_Def_op_a_xyk}
(\op[a]f)(x)& := \ossumint e^{i\bk\cdot(x-y)}a(x,y,\bk)f(y)\,\dbar y\\
& := \lim\limits_{\varepsilon\searrow 0}\sumint e^{i\bk\cdot(x-y)} \chi(\varepsilon \bk) a(x,y,\bk)f(y)\,\dbar y.
\end{split}
 \end{align}
\end{definition}

\begin{remark}
Let $N\in \N_0$ with $m+n<N\leq r$. If $r>0$, we can take $N\in\N$. In this case, integration by parts and the theorem of dominated convergence, applied to the sum-integral in Definition \ref{Def_op_a_xyk}, yield
\begin{align}
 (\op[a]f)(x)& := \lim\limits_{\varepsilon\searrow 0}\sumint e^{i\bk\cdot(x-y)} \chi(\varepsilon \bk) a(x,y,\bk)f(y)\,\dbar y \nonumber\\
   & = \lim\limits_{\varepsilon\searrow 0}\sumint \big(L_{\bk,y}^N e^{i\bk\cdot(x-y)}\big) \chi(\varepsilon \bk) a(x,y,\bk)f(y)\,\dbar y\nonumber\\
   & = \lim\limits_{\varepsilon\searrow 0}\sumint  e^{i\bk\cdot(x-y)}\chi(\varepsilon \bk) (L_{\bk,y}^t)^N\big( a(x,y,\bk)f(y)\big)\,\dbar y\nonumber\\
   & = \sumint  e^{i\bk\cdot(x-y)}(L_{\bk,y}^t)^N\big( a(x,y,\bk)f(y)\big)\,\dbar y. \label{Eq_reg_osc_sum_int}
  \end{align}
  Indeed, using \eqref{Eq-L_kytN} and the Leibniz rule, we obtain
  \begin{align*}
     (L_{\bk,y}^t)^N\big( a(x,y,\bk)f(y)\big) & = \sum\limits_{\substack{\gamma\in \N_0^n\\   \vert\gamma\vert\leq N}} c_{\gamma,N}\partial_y^\gamma\big( a(x,y,\bk)f(y)\big)\\
     & = \sum\limits_{\substack{\gamma\in \N_0^n\\   \vert\gamma\vert\leq N}} c_{\gamma,N}\sum\limits_{\substack{\beta\in \N_0^n\\   \beta\leq\gamma}}\binom{\gamma}{\beta}(\partial_y^\beta a(x,y,\bk))\partial_y^{\gamma-\beta}f(y).
  \end{align*}
Therefore, we have
  \begin{align*}
      \big\Vert e^{i\bk\cdot(x-y)}\chi(\varepsilon \bk) (L_{\bk,y}^t)^N\big( a(x,y,\bk)f(y)\big) \big\Vert & \leq \tilde{C}_N \Vert a\Vert_{m}^{(\rho,r)}\Vert f\Vert_{C^N(\T^n,E)}\langle\bk\rangle^{m-N}.
  \end{align*}
  The last inequality holds also true in the case $r=0$ with $N=0$. In this case $m<-n$.\
  Note that $N$ cannot be greater than $r$ due to the regularity of the symbol in $\T^{2n}$. The condition $m+n<N$ allows us to use dominated convergence and guarantees that the last sum-integral in \eqref{Eq_reg_osc_sum_int} is absolutely convergent. In order to choose $N\in \N$ such that $m+n<N\leq r$, we must have $r \geq \lfloor m  \rfloor + n + 1$. This is why we assume this condition in Definition \ref{Def_op_a_xyk}. Note also that the equation \eqref{Eq_reg_osc_sum_int} shows that Definition \ref{Def_op_a_xyk} is independent of the choice of the function $\chi$ and that   \eqref{Eq_reg_osc_sum_int} can be regarded as a \emph{regularization} of the oscillatory sum-integral.
\end{remark}

Throughout this section, operators associated with single symbols
$a(x,\bk)$ will occasionally be denoted by $a(X,D)$. Likewise, operators
associated with double symbols $a(x,y,\bk)$ will be denoted by
$a(X,Y,D)$.

If a double symbol is independent of one of its spatial variables,
we shall omit the corresponding formal variable. For instance, if
$a(x,y,\bk)=a(y,\bk)$ is independent of $x$, then the associated operator
will be denoted by $a(Y,D)$. Similarly, if $a(x,y,\bk)=a(x,\bk)$ is
independent of $y$, then we write $a(X,D)$. In particular, every ordinary symbol $a(x,\bk)$ may be regarded as the double symbol $a(x,y,\bk):=a(x,\bk)$, which is independent of $y$. Taking these conventions into account, the following result is true.

\begin{lemma} \label{2.3}
  Let $m_{1},\,m_2\in\R$, $\rho\in\N$, $r\in\N$ with $r\geq \max\big\{\lfloor m_2\rfloor +n + 1,\  \lfloor m_1+m_2\rfloor +n + 1\big\}$,  $a_1\in S^{m_1,\rho,r}(\T^n\times\Z^n,L(E))$ and $a_2\in S^{m_2,\rho,r}(\T^{2n}\times\Z^n,L(E))$ with $a_2$ independent of $x$. Let $a(x,y,\bk):=a_1(x,\bk)a_2(y,\bk)$ for $x,y\in\T^n$ and $\bk\in\Z^n$. Then, $a\in S^{m_1+m_2,\rho,r}(\T^{2n}\times\Z^n,L(E))$ and $a(X,Y,D)f=a_1(X,D)a_2(Y,D)f$ for $f\in C^\infty (\T^n,E)$.\ Here, $a_1(X,D)$ is evaluated in $a_2(Y,D)f$ by its Fourier series,
which converges absolutely and uniformly.
\end{lemma}
\begin{proof}
  By straightforward calculations, using the discrete Leibniz formula \eqref{4-0a}, we have that $a\in S^{m_1+m_2,\rho,r}(\T^{2n}\times\Z^n,L(E))$. Now, let $f\in C^\infty(\T^n,E)$ and take $N\in\N$ with
\[\max\big\{m_1+m_2+n,\ m_2+n \big\}<N\leq r.\]
Note that by \eqref{Eq_reg_osc_sum_int}, for $\eta\in\T^n$ we have
\begin{equation*}
    (a_2(Y,D)f)(\eta) = \sum\limits_{\tilde{\bk}\in\Z^n}\int_{\T^n}e^{i\tilde{\bk}\cdot(\eta-y)}{(L_{\tilde{\bk},y}^t)^N}[a_2(y,\tilde{\bk})f(y)]\,\dbar y.
\end{equation*}
Put $u:=a_2(Y,D)f$ and
\[
 c_{\tilde{\bk}}:=
 \int_{\T^n}e^{-i\tilde{\bk}\cdot y}
    (L_{\tilde{\bk},y}^t)^N[a_2(y,\tilde{\bk})f(y)]\,\dbar y.
\]
Therefore,
$$u(\eta)=\sum\limits_{\tilde{\bk}\in\Z^n}e^{i\tilde{\bk}\cdot \eta}c_{\tilde{\bk}}\qquad (\eta\in\T^n).$$
Since
$$\Vert c_{\tilde{\bk}}\Vert\le \tilde{C}_N \Vert a_2\Vert_{m_2}^{(\rho,r)}\Vert f\Vert_{C^N(\T^n,E)}\langle\tilde\bk\rangle^{m_2-N} \qquad (\tilde\bk\in\Z^n),$$
the previous Fourier series converges absolutely and uniformly, so that $\widehat u(\bk)=c_{\bk}$ for $\bk\in\Z^n$.
Moreover,
\[
 \sum_{\bk\in\Z^n}\sup_{x\in\T^n}\Vert a_1(x,\bk)c_{\bk}\Vert
 \le C_f\sum_{\bk\in\Z^n}
              \langle\bk\rangle^{m_1+m_2-N}<\infty.
\]
Thus, the Fourier series defining $a_1(X,D)u$ also converges absolutely
and uniformly, and we obtain
\begin{align*}
 \big[a_1(X,D)a_2(Y,D)f\big](x)&= [a_1(X,D)u](x)=\sum_{\bk\in\Z^n}e^{i\bk\cdot x}a_1(x,\bk)\hat{u}(\bk)\\
 &=\sum_{\bk\in\Z^n}e^{i\bk\cdot x}a_1(x,\bk)c_{\bk}\\
 & = \sum\limits_{\bk\in\Z^n}\int_{\T^n}e^{i\bk\cdot(x-y)}{(L_{\bk,y}^t)^N}[a_1(x,\bk)a_2(y,\bk)f(y)]\dbar y\\
 & = \sum\limits_{\bk\in\Z^n}\int_{\T^n}e^{i\bk\cdot(x-y)}{(L_{\bk,y}^t)^N}[a(x,y,\bk)f(y)]\dbar y\\
 & = \big(a(X,Y,D)f\big)(x).
\end{align*}
\end{proof}

\begin{proposition}\label{Prop. CM-Ck-continuity} Let $m\in\mathbb R$, $\rho\in\mathbb N_0$, $r,k,M\in\mathbb N_0$, and let $a\in S^{m,\rho,r}(\mathbb T^n\times\mathbb Z^n,L(E))$. Assume that
$r\ge k$ and $M>n+m+k$ . Then the pseudodifferential operator
\[
\operatorname{op}[a]: C^\infty(\mathbb T^n,E) \to C^r (\mathbb T^n,E)
\]
extends uniquely to a bounded linear operator
\[
\operatorname{op}[a]: C^M(\mathbb T^n,E) \to C^k(\mathbb T^n,E).
\]
More precisely, there exists a constant \(C>0\) such that
\[
\|\operatorname{op}[a]f\|_{C^k(\T^n,E)} \le C\|f\|_{C^M(\T^n,E)} \qquad \text{for all} \,  f\in C^M(\mathbb T^n,E).
\]
\end{proposition}

\begin{proof}
Let $f\in C^M(\mathbb T^n,E)$. We first estimate its Fourier coefficients. Let $L_{\bk,y}$ be given by \eqref{Eq-L_ky}. For
\({\bk}\in\mathbb Z^n\), we have
\begin{align*}
\widehat f(\bk)
&=
\int_{\mathbb T^n}
e^{-i\bk\cdot y}
f(y)\,\dbar y=
\int_{\mathbb T^n}
e^{-i\bk\cdot y}
\bigl(L_{\bk,y}^{\,t}\bigr)^M
f(y)\,\dbar y,
\end{align*}
\noindent where $L_{\bk,y}^t$ is given by \eqref{Eq-L_kyt}. Taking into account \eqref{Eq-L_kytN},  we obtain
\begin{equation}\label{eq:Fourier-CM-estimate}
\|\widehat f(\bk)\|
\le
C_M
\langle\bk\rangle^{-M}
\|f\|_{C^M(\T^n,E)}.
\end{equation}
Consider now the Fourier series
\[
\operatorname{op}[a]f(x)
=
\sum_{\bk\in\mathbb Z^n}
e^{i\bk\cdot x}
a(x,\bk)
\widehat f(\bk).
\]
Let \(\alpha\in\mathbb N_0^n\) with $|\alpha|\le k$. By the Leibniz rule,
\begin{equation} \label{eq:derivative-op-a}
\partial_x^\alpha
\Bigl(
e^{i\bk\cdot x}
a(x,\bk)
\widehat f(\bk)
\Bigr) =\sum_{\beta\le\alpha}
\binom{\alpha}{\beta}
(i\bk)^{\alpha-\beta}
e^{i\bk\cdot x}
\partial_x^\beta a(x,\bk)
\widehat f(\bk).
\end{equation}
Since \(a\in S^{m,\rho,r}\) and \(r\ge k\), for every
\(|\beta|\le k\) we have
\[
\sup_{x\in\mathbb T^n}
\|
\partial_x^\beta a(x,\bk)
\|_{L(E)}
\le
C_\beta
\langle\bk\rangle^m.
\]
Combining this estimate with
\eqref{eq:Fourier-CM-estimate}, we obtain
\begin{align*}
\sup_{x\in\mathbb T^n}
\left\|
\partial_x^\alpha
\Bigl(
e^{i\bk\cdot x}
a(x,\bk)
\widehat f(\bk)
\Bigr)
\right\|_E
&\le
C_\alpha
\sum_{\beta\le\alpha}
\langle\bk\rangle^{|\alpha-\beta|}
\langle\bk\rangle^m
\langle\bk\rangle^{-M}
\|f\|_{C^M(\T^n,E)}
\\
&\le
C_\alpha
\langle\bk\rangle^{m+|\alpha|-M}
\|f\|_{C^M(\T^n,E)}
\\
&\le
C_\alpha
\langle\bk\rangle^{m+k-M}
\|f\|_{C^M(\T^n,E)}.
\end{align*}

The assumption $M>n+m+k$ is equivalent to $m+k-M<-n$. Therefore,
\[
\sum_{\bk\in\mathbb Z^n}
\langle\bk\rangle^{m+k-M}
<\infty.
\]
It follows from the Weierstrass \(M\)-test that, for every
\(|\alpha|\le k\), the series
\[
\sum_{\bk\in\mathbb Z^n}
\partial_x^\alpha
\Bigl(
e^{i\bk\cdot x}
a(x,\bk)
\widehat f(\bk)
\Bigr)
\]
converges absolutely and uniformly with respect to
\(x\in\mathbb T^n\).

Hence, the Fourier series defining
\(\operatorname{op}[a]f\) belongs to
\(C^k(\mathbb T^n,E)\), and termwise differentiation is valid up to
order \(k\). Moreover,
\begin{align*}
\|\operatorname{op}[a]f\|_{C^k(\T^n,E)}
&=
\max_{|\alpha|\le k}
\sup_{x\in\mathbb T^n}
\left\|
\partial_x^\alpha
\operatorname{op}[a]f(x)
\right\|_E
\\
&\le
C
\sum_{\bk\in\mathbb Z^n}
\langle\bk\rangle^{m+k-M}
\|f\|_{C^M(\T^n,E)}
\\
&\le
C'\|f\|_{C^M(\T^n,E)}.
\end{align*}
Thus,
\[
\operatorname{op}[a]:
C^M(\mathbb T^n,E)
\to
C^k(\mathbb T^n,E)
\]
is bounded.

Since \(C^\infty(\mathbb T^n,E)\) is dense in
\(C^M(\mathbb T^n,E)\), this bounded operator is the unique continuous
extension of the original operator defined on
\(C^\infty(\mathbb T^n,E)\).
\end{proof}

\section{Generation of semigroups for parabolic pseudodifferential operators}

Now we will investigate the generation of analytic semigroups by realizations of toroidal pseudodifferential operators in Sobolev and Besov spaces. For this, one of the key ingredients is an estimate of the inverse symbol $(a(x,\bk)+\lambda)^{-1}$ for a parabolic symbol $a\in S^{m,\rho,r}$.\ We therefore start with a remark on the discrete derivatives.

\begin{remark}
  \label{4.0}
  Let $a\colon \T^n\times\Z^n\to L(E)$ and $\lambda\in\C$ such that $a(x,\bk)+\lambda$ is invertible for all $\bk\in\Z^n$ and $x\in\T^n$. Then, the discrete difference $\Delta_\bk^\gamma (a(x,\bk)+\lambda)^{-1}$, with $\vert \gamma\vert>0$, can be written as a finite linear combination of terms of the form
  \begin{equation}
   \label{4-0}
   \big(a(x, \bk+\beta_0^{(j)})+\lambda\big)^{-1} \prod_{l=1}^j \Big[ (\Delta_\bk^{\alpha^{(l)}} a)(x,\bk+\beta_1^{(l)}) \big( a(x,\bk+\beta_2^{(l)})+\lambda\big)^{-1}\Big],
   \end{equation}
  where $\alpha^{(l)}\in\N_0^n\smallsetminus\{0\}$, {$\beta_i^{(l)}\in\N_0^n$}, $j\in\{1,\dots,|\gamma|\}$, and  $\sum_{l=1}^j |\alpha^{(l)}| = |\gamma|$.

  This statement follows by induction in a straightforward way, based on the discrete Leibniz formula \eqref{4-0a}.
  \end{remark}

Throughout this section, we fix $m\in (0,\infty)$, $r\in\N_0$ and $\rho\in\N$.

\begin{definition}
  \label{4.1}
  The symbol $a\in S^{m,\rho,r}$ is called parabolic with constants $\omega\ge 0$ and $\kappa>0$ if for all $(x,\bk,\lambda)\in\T^n\times \Z^n\times \C$ with $\Re\lambda\ge 0$  and $\big|\big(\bk,\vert\lambda\vert^{1/m}\big)\big|\ge \omega$ we have that $a(x,\bk)+\lambda\colon E\to E$ is bijective and
  \[ \big\| \big(a(x,\bk)+ {\lambda}\big)^{-1} \big\|_{L(E)} \le \kappa \big\langle\bk, \vert\lambda\vert^{1/m}\big\rangle^{-m}.\]
  Here, $\langle \bk,\mu\rangle := (1+|\bk|^2+\mu^2)^{1/2}$ and $|(\bk,\mu)| := (|\bk|^2+\mu^2)^{1/2}$ for $\mu\ge 0$.
  We write $P_{\omega,\kappa}S^{m,\rho,r}$ for the set of all symbols in $S^{m,\rho,r}$, which are parabolic with constants $\omega$ and $\kappa$ and endow this set with the topology induced by the topology in $S^{m,\rho,r}$.
\end{definition}

\begin{definition}
 For $a\in P_{\omega,\kappa}S^{m,\rho,r}$ and {$\lambda\in\C$ with $\Re\lambda\ge0$} let $a_\lambda$ be {defined} by $a_\lambda(x,\bk):=\lambda + a(x,\bk)$ for $(x,\bk)\in\T^n\times \Z^n$. {Furthermore, let}  $\psi\in C^\infty(\R^{n+1})$ with $0\le \psi\le 1$, $\psi(\xi,\mu)=0$ for ${|(\xi,\mu)|}\le \frac12$ and $\psi(\xi,\mu)=1$ for ${|(\xi,\mu)|}\ge 1$ be a fixed {cut-off function}, and $\omega_0>\omega$. For $\mu:=\vert \lambda\vert^{1/m}$ we set
\[ a_\lambda^\#(x,\bk) :=  \begin{cases}
  \psi\big(\frac\bk{2\omega_0},\frac\mu{2\omega_0}\big) a_\lambda(x,\bk)^{-1}, & \text{if } x\in\T^n,\, |\bk,\mu|\ge \omega_0,\\
  0, & \text{if }x\in\T^n, \,|\bk,\mu|\le \omega_0.
\end{cases}\]
\end{definition}

\begin{remark}
 Note that in the definition above $a_\lambda^\#(x,\bk)=0$ if $|\bk,\mu|\le \omega_0$ and $a_\lambda^\#(x,\bk)=a_\lambda(x,\bk)^{-1}$ if $|\bk,\mu|\ge 2\omega_0$ for $x\in\T^n$. In particular, we have
 \[a_\lambda^\#(x,\bk)=a_\lambda(x,\bk)^{-1} \qquad (x\in\T^n, \bk\in\Z^n, \mu\ge 2\omega_0).\]
\end{remark}

\begin{lemma}
  \label{Lemma2.5}
  Let $\mathcal A\subset P_{\omega,\kappa}S^{m,\rho,r}$ be bounded, i.e. there exists a constant $K>0$ such that $\Vert a\Vert_{m}^{(\rho,r)}\le K\;(a\in\mathcal A)$. Then there exists a constant $C=C(K,\kappa,\rho,r)>0$ {(with $n,m,\omega_0$ and $\psi$ fixed)} such that, with $\mu=\vert \lambda\vert^{1/m}$, the following holds:
  \[ \| \partial_x^\beta a_\lambda^\#(x,\bk)\|_{L(E)} \le C\langle\bk,\mu\rangle^{-m},\]
  \[ \| \Delta_\bk^\alpha\partial_x^\beta a_\lambda^\#(x,\bk)\|_{L(E)} \le C\langle\bk,\mu\rangle^{-2m}\langle\bk\rangle^{m-\vert\alpha\vert},\]
  for all $x\in\T^n, (\bk,\mu)\in\Z^n\times {[0,\infty)},0<|\alpha|\le \rho,\,  |\beta|\le r$, and $a\in\mathcal A$, {where $\lambda\in\C$ satisfies $\Re\lambda\ge0$ and $|\lambda|=\mu^m$}, which implies that $\{a_\lambda^\#:a\in\mathcal A, \lambda\in\C,\ \Re\lambda\ge0\}\subset S^{-m,\rho,r}$ is bounded.
\end{lemma}

\begin{proof}
Put $R:=2\omega_0+\rho+1$. On $|(\bk,\mu)|\le R$, all arguments
$\bk+\eta$ occurring in the differences satisfy $|\eta|\le\rho$.
Where the cut-off is nonzero, parabolicity bounds the inverse uniformly by iterating
the identity $\partial_{x_j}a_\lambda^{-1}
=-a_\lambda^{-1}(\partial_{x_j}a)a_\lambda^{-1}$, which gives
uniform bounds for the spatial derivatives. Together with the fixed
cut-off and the uniform bound for the symbol, this proves both asserted
estimates in the bounded region where $|(\bk,\mu)|\leq R$, uniformly for $a\in\mathcal A$.
For $|(\bk,\mu)|>R$, all shifted cut-offs equal one.
Let $0<|\alpha|\le\rho$ and $|\beta|\le r$. Then,
$\Delta_\bk^\alpha\partial_x^\beta a_\lambda(x,\bk)^{-1}$ is a finite
linear combination of terms of the form
  \begin{equation}
    \label{3-1}
        a_\lambda\big(x,\bk+\eta_0^{(p)}\big)^{-1} {\prod_{i=1}^p \Big[\big(\Delta_\bk^{\alpha^{(i)}}\partial_x^{\beta^{(i)}} a_\lambda\big)\big(x,\bk+\eta_1^{(i)}\big)a_\lambda\big(x,\bk+\eta_2^{(i)}\big)^{-1}\Big]}
       \end{equation}
  with $1\le p \le |\alpha|+|\beta|$ and with {$\alpha^{(i)},\beta^{(i)}\in\N_0^n$, $|\alpha^{(i)}|+|\beta^{(i)}|>0$} $\;(i=1,\dots,p)$ satisfying $\sum_{i=1}^p \alpha^{(i)}=\alpha$, $\sum_{i=1}^p \beta^{(i)}=\beta$, and $\eta_j^{(i)}$ are suitable multi-indices {of length at most $\rho$}. Using Peetre's inequality, the norm of the operator in \eqref{3-1} can be estimated by
   \begin{align*}
  &\begin{multlined}[t][0.9\linewidth]
  \big\Vert  a_\lambda\big(x,\bk+\eta_0^{(p)}\big)^{-1}\big\Vert_{L(E)} \prod_{i=1}^p \big\Vert\big[\Delta_\bk^{\alpha^{(i)}}\partial_x^{\beta^{(i)}} a_\lambda\big(x,\bk+\eta_1^{(i)}\big)\big]\big\Vert_{L(E)}\\
     \cdot \big\Vert a_\lambda\big(x,\bk+\eta_2^{(i)}\big)^{-1}\big\Vert_{L(E)}
  \end{multlined}\\
    & \qquad\le C_p\kappa^{p+1}\langle\bk,\mu\rangle^{-m(p+1)}\prod_{i=1}^p \Vert a\Vert_m^{(\rho,r)}\langle\bk\rangle^{m-|\alpha^{(i)}|}\\
    & \qquad= C_p \kappa^{p+1}\big(\Vert a\Vert_m^{(\rho,r)}\big)^p\langle \bk,\mu\rangle^{-mp-m} \langle\bk\rangle^{pm-|\alpha|}\\
    & \qquad= C_p \kappa^{p+1}\big(\Vert a\Vert_m^{(\rho,r)}\big)^p\langle \bk,\mu\rangle^{-2m}\langle \bk,\mu\rangle^{-m(p-1)}\langle\bk\rangle^{pm-|\alpha|}\\
    & \qquad\le C_p \kappa^{p+1}\big(\Vert a\Vert_m^{(\rho,r)}\big)^p\langle \bk,\mu\rangle^{-2m}\langle\bk\rangle^{m-\vert \alpha\vert}.
  \end{align*}
  Summing up over all terms of the form \eqref{3-1}, we obtain
  \[ \| \Delta_\bk^\alpha\partial_x^\beta a_\lambda(x,\bk)^{-1}\|_{L(E)} \le C\langle \bk,\mu\rangle^{-2m}\langle\bk\rangle^{m-\vert \alpha\vert}\]
  with a constant $C=C(K,\kappa,\rho,r)$.\\

 Now, if $\alpha =\beta = 0$, then
  \[ \big\| \Delta_\bk^\alpha\partial_x^\beta  a_\lambda(x,\bk)^{-1} \big\|_{L(E)} = \big\|  a_\lambda(x,\bk)^{-1} \big\|_{L(E)} \le \kappa\langle\bk,\mu\rangle^{-m}.\]
  If $\alpha=0$ and $\beta\in \N_0^n\smallsetminus\{0\}$, then $\Delta_\bk^\alpha\partial_x^\beta a_\lambda(x,\bk)^{-1} = \partial_x^\beta a_\lambda(x,\bk)^{-1}$ is a finite linear combination of {terms} of the form
  \[ a_\lambda(x,\bk)^{-1}\prod\limits_{i=1}^p \Big[\big(\partial_x^{\beta^{(i)}}a(x,\bk)\big)a_\lambda(x,\bk)^{-1}\Big]\]
  with $1\le p\le \vert\beta\vert$ and $\beta^{(i)}\in\N_0^n\smallsetminus\{0\}\ (i=1,\dots,p)$ satisfying $\sum\limits_{i=1}^p \beta^{(i)} = \beta$. Then,
  \begin{align*}
   \big\| \Delta_\bk^0\partial_x^\beta a_\lambda(x,\bk)^{-1}\big\|_{L(E)} & \le {C_\beta\sum_{p=1}^{|\beta|}}\kappa\langle \bk,\mu\rangle^{-m}\prod_{i=1}^p\big( \Vert a\Vert_m^{{(\rho,r)}}\langle \bk\rangle^m\big)\big(\kappa\langle \bk,\mu\rangle^{-m}\big)\\
   & \le C_\kappa {\sum_{p=1}^{|\beta|}}\big(\Vert a\Vert_m^{(\rho,r)}\big)^p\langle \bk,\mu\rangle^{-m}\prod\limits_{i=1}^p \underbrace{\langle\bk\rangle^m\langle\bk,\mu\rangle^{-m}}_{\le 1}\\
   & \le C\langle \bk,\mu \rangle^{-m}.
  \end{align*}

\end{proof}

\begin{remark}\label{rem:sector-extension}
There are $\delta>0$ and $\tilde{R}>0$, depending only on the
parabolicity constants and $m$, such that the previous
estimates for the inverse symbol also hold, without a cut-off, for
$|\lambda|\ge \tilde{R}$ and
$|\arg\lambda|\le\pi/2+\arctan\delta$.
Indeed, only the case $\Re\lambda<0$ needs to be considered. Set
$\lambda_0=i\Im\lambda$. Then
$|\lambda-\lambda_0|\le\delta|\lambda_0|$ and, for $\tilde{R}$ large,
parabolicity gives
\[
\big\| (\lambda-\lambda_0)\big(a(x,\bk)+\lambda_0\big)^{-1}\big\|
\le\delta\kappa |\lambda_0|
\big\langle\bk,|\lambda_0|^{1/m}\big\rangle^{-m}
\le\delta\kappa<\tfrac12.
\]
A Neumann series and the relation between $|\lambda_0|$ and $|\lambda|$
give the inverse estimate with $\mu=|\lambda|^{1/m}$. The same formulas as above then yield the estimates for the spatial and
discrete derivatives, uniformly in the enlarged sector.
\end{remark}

\begin{lemma}\label{Lemma 2.6}\
 \begin{enumerate}
 \item[\textnormal{a)}]Let $\mu_0>0$ and $\varphi\in\mathscr S(\R)$ with $\varphi(0)=1$. For $\epsilon\in (0,1)$, set $\varphi_\epsilon(x) := \varphi(\epsilon x)\;(x\in\R)$ and $\chi_\epsilon(\bk, \mu) := \varphi_\epsilon \big((|\bk|^2+ \mu^2)^{1/2}\big)\;(\bk\in\Z^n,\, \mu>0)$. Then, for all $\gamma\in\N_0^n$  we have
  \[  {\big| \Delta_{\bk} ^\gamma \chi_\epsilon(\bk,\mu)\big| \le C \langle\bk, \mu\rangle^{-|\gamma|}\quad (\bk\in \Z^n, \mu\ge \mu_0)}\]
  with a constant $C$ depending on $\gamma$ and $\mu_0$ {and on $\varphi$} but not on $\epsilon$ or $\mu$.
  \item[\textnormal{b)}] Let $\mu_0>0$ and $\big\{b_\mu : \mu\ge\mu_0\big\} \subset S^{-m,\rho,r}(\T^{2n}\times \Z^n, L(E))$ bounded with
  \[ {\langle \bk,\mu\rangle^m\Vert \partial_{x,y}^\beta b_\mu(x,y,\bk)\Vert_{L(E)} + \langle \bk,\mu\rangle^{2m}\langle \bk\rangle^{\vert \gamma\vert - m}\Vert \Delta_\bk^\gamma\partial_{x,y}^\beta b_\mu(x,y,\bk)\Vert_{L(E)} \le C}\]
  for all $x,y \in\T^n$, $(\bk,\mu)\in \Z^n\times [\mu_0, \infty)$ and $\gamma\in\N_0^n$ with $0<\vert \gamma\vert\le \rho$ {and $|\beta|\le r$}. Then
\begin{equation*}
\Big\Vert \Delta_\bk^\gamma\big[ \chi_\epsilon(\bk,\mu){\partial_{x,y}^\beta}b_\mu(x,y,\bk)\big]\Big\Vert_{L(E)}
\le C_1 \langle \bk,\mu\rangle^{-m-\vert\gamma\vert} + C_2\langle \bk,\mu\rangle^{-2m}\langle \bk\rangle^{m-\vert\gamma\vert}
\end{equation*}
for all $x,y\in\T^n$, $(\bk,\mu)\in\Z^n\times[\mu_0,\infty)$ and $\vert\gamma\vert \le \rho$ {and $|\beta|\le r$}.
\end{enumerate}
\end{lemma}
\begin{proof}
 $a)$ Let $g(\xi, \mu):= (|\xi|^2+ \mu^2)^{1/2}$ for $(\xi,\mu)\in \R^n\times(0,\infty)$ so that $\chi_\epsilon(\bk,\mu) = (\varphi_\epsilon\circ g)(\bk,\mu)$. {Since $\varphi_\epsilon^{(k)}(x)=\epsilon^k\varphi^{(k)}(\epsilon x)$ and $\epsilon\langle x\rangle\le\langle\epsilon x\rangle$, for the $k$-th derivative of $\varphi_\epsilon$}, we have
\begin{equation}
  \label{4-8}
  {|\varphi_\epsilon^{(k)}(x)|}\le c_k \langle x\rangle^{-k}\quad (k\in\N_0,\, x\in\R)
\end{equation}
with a constant $c_k$ independent of $\epsilon\in (0,1)$. On the other hand, $g$ is positively homogeneous of degree $1$ in $(\xi, \mu)$, i.e. $g(t\xi, t\mu) = t g(\xi, \mu)$ holds for $(\xi,\mu)\in\R^n\times(0,\infty)$  and $t>0$. Therefore, the $\alpha$-th derivative is positively homogeneous of degree $1-|\alpha|$, i.e.,
\[ \partial^\alpha_{ {\xi}} g(t\xi, t\mu) = t^{1-|\alpha|} (\partial^\alpha_{ {\xi}} g)(\xi, \mu)\qquad  \big((\xi,\mu)\in\R^n\times(0,\infty),\, t>0\big).\]
{As $g$ extends smoothly to $\R^{n+1}\setminus\{0\}$, all its derivatives are bounded on the unit sphere}.\ Setting $t:= (|\xi|^2+ \mu^2)^{1/2}$, this implies that
\begin{equation}
  \label{4-9}
  |\partial^\alpha_{ {\xi}} g(\xi, \mu)| \le C_\alpha (|\xi|^2+ \mu^2)^{(1-|\alpha|)/2}\qquad  { \big((\xi,\mu)\in\R^n\times(0,\infty)\big).}
\end{equation}
To estimate $\partial^\gamma_{ {\xi}}\chi_\epsilon$ {for $|\gamma|>0$}, we apply the generalized chain rule which can be formulated as follows (by repeated application of the ordinary product and chain rules, see also \cite{bahouri-chemin-danchin11} for more details):\ The $\gamma$-th derivative $[\partial^\gamma_{ {\xi}} (\varphi_\epsilon\circ g)](\xi,\mu)$ is a finite linear combination of terms of the form
\[ \varphi_\epsilon^{(k)}\big(g(\xi,\mu)\big)\cdot \big(\partial^{\alpha^{(1)}}_{ {\xi}} g\big)(\xi,\mu)\cdot\ldots\cdot  \big(\partial^{\alpha^{(k)}}_{ {\xi}} g\big)(\xi,\mu), \]
where $k\in \{1,\dots,|\gamma|\}$ and $\alpha^{(1)},\dots, \alpha^{(k)}\in{\N_0^n\setminus\{0\}}$ with $|\alpha^{(1)}|+\ldots+|\alpha^{(k)}| = |\gamma|$.

From this and the estimates \eqref{4-8} and \eqref{4-9}, we obtain
\[ |\partial^\gamma_{ {\xi}} \chi_\epsilon(\xi, \mu)|\le C_\gamma (|\xi|^2+ \mu^2)^{-|\gamma|/2}.\]
Here the constant $C_\gamma$ depends on $\gamma$ and ${\varphi}$ but not on $\epsilon$, $\xi$ or $\mu$.

To pass to discrete differences, we consider first $d\coloneqq|\gamma|>0$. Note that we can write
$\gamma=\bde_{j_1}+\cdots+\bde_{j_d}$, where $j_1,\dots,j_d\in\{1,\dots,n\}$.
Repeated use of the fundamental theorem of calculus gives
\[
\Delta_\bk^\gamma\chi_\epsilon(\bk,\mu)
=\int_{[0,1]^d}(\partial_\xi^\gamma\chi_\epsilon)
  \Big(\bk+\sum_{\nu=1}^d t_\nu\bde_{j_\nu},\mu\Big)\,dt_1\cdots dt_d.
\]
For $\mu\ge\mu_0$, the previous derivative estimate and Peetre's
inequality bound the integrand by $C_{\gamma,\mu_0}
\langle\bk,\mu\rangle^{-d}$. This proves the assertion for $d>0$.\
The case $d=0$ follows directly from $\|\varphi\|_\infty<\infty$.

$b)$ The estimate follows from $a)$ by {the discrete Leibniz formula} \eqref{4-0a}.
\end{proof}

It was shown in \cite{barraza-denk-hernandez-nau} that for $m>0$ and $\mu\ge\mu_0>0$ we have

\begin{align}
    \sum_{\substack{\ell\in\mu^{-1}\Z^n\setminus \{0\} \\ |\ell|_\infty\le 1}} |\ell|^{\theta-n} \mu^{-n} & \le C_\theta && (\theta>0) ,\label{2-6}\\
    \sum_{\substack{\ell\in\mu^{-1}\Z^n \\ |\ell|_\infty> 1}} |\ell|^{\theta-n-m}\mu^{-n} & \le C_\theta && \big(\theta\in (0,m)\big) \label{2-7}.
  \end{align}
These estimates will be useful for the proof of the following lemma.
\begin{lemma}\label{Lemma2.7}
Let $\theta_0:=\frac{1}{2}\min\{1, m\}$, and $\theta_1:=1-\theta_0$.   Then, for all $\mu_0>0$ there exists a constant $C>0$ such that for all $\mu\geq \mu_0$ and $i=0, 1$, the following estimate holds:
 \[ \sum\limits_{\bk\in\Z^n} \vert \bk\vert^{\theta_i}\langle \bk,\mu\rangle^{-m-(n+i)} \le C\mu^{\theta_i}\mu^{-(m+i)}. \]
 \end{lemma}
\begin{proof}
 At first note that, for $\theta>0$,
 \[ \sum\limits_{\bk\in\Z^n} \vert \bk\vert^{\theta}\langle \bk,\mu\rangle^{-m-n} \le C\mu^{\theta}\mu^{-m} \]
 is equivalent to
 \[ \sum\limits_{\ell\in\mu^{-1}\Z^n}\vert \ell\vert^\theta\big(1+\mu^2+\vert\mu\ell\vert^2\big)^{\frac{-m-n}{2}} \le C\mu^{-m}.\]
 Now, due to \eqref{2-6}, for $\theta>0$ we have
 \begin{align*}
  \sum_{\substack{\ell\in\mu^{-1}\Z^n\setminus \{0\} \\ |\ell|_\infty\le 1}} & \mu^m  \vert \ell \vert^\theta\big(1+\mu^2+\vert\mu\ell\vert^2\big)^{\frac{-m-n}{2}}\\
  & \le \sum_{\substack{\ell\in\mu^{-1}\Z^n\setminus \{0\} \\ |\ell|_\infty\le 1}}\underbrace{\mu^m\big(1+\mu^2+\vert\mu\ell\vert^2\big)^{\frac{-m}{2}}}_{\le 1}\vert \ell\vert^\theta \underbrace{\big(1+\mu^2+\vert\mu\ell\vert^2\big)^{\frac{-n}{2}}}_{\le \mu^{-n}\vert\ell\vert^{-n}}\\
  & \le \sum_{\substack{\ell\in\mu^{-1}\Z^n\setminus \{0\} \\ |\ell|_\infty\le 1}} \vert \ell\vert^{\theta-n}\mu^{-n} \le C.
 \end{align*}
{On the other hand}, due to \eqref{2-7}, for $0<\theta<m$ we also have
\begin{equation*}
 \sum_{\substack{\ell\in\mu^{-1}\Z^n \\ |\ell|_\infty > 1}} \mu^{m}\vert \ell\vert^\theta \underbrace{\big(1+\mu^2+\vert\mu\ell\vert^2\big)^{\frac{-m-n}{2}}}_{\le \mu^{-m-n}\vert\ell\vert^{-m-n}} \le \sum_{\substack{\ell\in\mu^{-1}\Z^n \\ |\ell|_\infty > 1}} \vert\ell\vert^{\theta-n-m}\mu^{-n} \le C.
\end{equation*}
Therefore, it follows for $0<\theta<m$ that
\begin{equation}\label{2-8}
 \sum\limits_{\ell\in\mu^{-1}\Z^n} \mu^m\vert \ell\vert^{\theta}\big(1+\mu^2+\vert\mu\ell\vert^2\big)^{\frac{-m-n}{2}} \le C.
\end{equation}
Then, setting $\theta=\theta_0$ in \eqref{2-8}, we obtain
\[ \sum\limits_{\bk\in\Z^n} \vert \bk\vert^{\theta_0}\langle \bk,\mu\rangle^{-m-n} \le C\mu^{\theta_0}\mu^{-m}. \]
Since $0<\theta_1<1<m+1$ and due to the fact that \eqref{2-8} holds for all $m>0$ and $0<\theta<m$, we obtain also
\begin{equation*}
 \sum\limits_{\ell\in\mu^{-1}\Z^n} \mu^{m+1}\vert \ell\vert^{\theta_1}\big(1+\mu^2+\vert\mu\ell\vert^2\big)^{\frac{-m-1-n}{2}} \le C,
\end{equation*}
which is equivalent to
\[ \sum\limits_{\bk\in\Z^n} \vert \bk\vert^{\theta_1}\langle \bk,\mu\rangle^{-m-(n+1)} \le C\mu^{\theta_1}\mu^{-(m+1)}. \]
\end{proof}

\begin{remark}
  \label{Remark 4.7 JFAA}
  In the proof of the following lemma, we will use the elementary inequality
  \begin{equation}
    \label{4-10}
    |\eta|^N \le C  \sum_{|\gamma|=N} \big| (e^{-i\eta}-1)^\gamma \big|\quad (\eta\in\T^n),
  \end{equation}
  for $N\in\N$, where we have set $(e^{-i\eta}-1)^\gamma := (e^{-i\eta_1} -1)^{\gamma_1}\cdot\ldots\cdot (e^{-i\eta_n}-1)^{\gamma_n}$.\\
To prove \eqref{4-10}, we consider the function $f(\eta) := \sum_{|\gamma|=N} | (e^{-i\eta}-1)^\gamma |$. Obviously, $f$ is a continuous function on $\T^n$ and has no zeros for $|\eta| \ge\frac\pi 2$ (choose $\gamma_j := N$ for $|\eta_j| = \max\{ |\eta_1| ,\dots,|\eta_n|\}$). Therefore $f(\eta)\ge C >0\;(|\eta|\ge\frac\pi 2)$ which implies \eqref{4-10} for $|\eta|\ge \frac\pi 2$.\\
For $|\eta|\le\frac\pi 2$, we use
  \[ |(e^{-i\eta}-1)^\gamma|  = \prod_{j=1}^n |e^{-i\eta_j}-1|^{\gamma_j} \ge \prod_{j=1}^n |\sin(\eta_j)|^{\gamma_j} \ge C\prod_{j=1}^n |\eta_j|^{\gamma_j} . \]
  Now \eqref{4-10} follows from
  \[   |\eta|^N \le C\sum_{|\gamma|=N} |\eta^\gamma| = C\sum_{|\gamma|=N} \prod_{j=1}^n |\eta_j|^{\gamma_j}.\]
\end{remark}
The following result is the key estimate in the proof of the generation of an analytic semigroup in Sobolev spaces.

\begin{lemma}
  \label{Lemma2.9}
  Let $\mu_0>0$, $\rho\ge n+1$,  $\big\{b_\mu\,:\,\mu\ge\mu_0 \big\}\subset S^{-m,\rho,r}(\T^{2n}\times\Z^n,L(E))$ and $\chi_\epsilon$ both as in Lemma \ref{Lemma 2.6}. For $0<\epsilon<1$, $(x,\eta)\in\T^{2n}$, and $\mu\ge\mu_0$, define
  \[ K_\epsilon (x,\eta,\mu) := \sum_{\bk\in\Z^n} e^{i\bk \cdot\eta} \chi_\epsilon(\bk, \mu)b_\mu(x,x-\eta,\bk)\]
  and
  \[ {g_\mu}(\eta):= \frac{  |\mu\eta| ^{\theta_0} + |\mu\eta| ^{\theta_1}}{ |\eta| ^n (1+|\mu\eta| )},\]
  where $\theta_0$ and $\theta_1$ are as in Lemma~\ref{Lemma2.7}. Then,
  \begin{equation}\label{4-11}
   \mu^m\, \|K_\epsilon(x,\eta,\mu)\|_{L(E)} \le C \, {g_\mu}(\eta) \quad (x\in\T^n,\ {\eta\in\T^n\setminus\{0\}},\, \mu\ge\mu_0,\, 0<\epsilon<1)
   \end{equation}
  with $C$ being independent of $\epsilon$, $\eta$, and $\mu$. Furthermore,
    \begin{equation}\label{4-12}
  \mu^m\,\|K_\epsilon(x,\cdot,\mu)\|_{L^1(\T^n,L(E))} \le C \quad (\mu\ge\mu_0,\,x\in\T^n,\,0<\epsilon<1).
  \end{equation}
  Moreover, {for each $\mu\ge\mu_0$, there exists a strongly measurable function $K(\cdot,\cdot,\mu)\colon \T^{2n}\to L(E)$} with $K_\epsilon(x,\eta,\mu)\to K(x,\eta,\mu)$  as $\epsilon \searrow 0$ pointwise almost everywhere, and \eqref{4-11} and \eqref{4-12} hold with $K_\epsilon$ being replaced by $K$.
\end{lemma}
\begin{proof}
   {Write $\bar\Delta_{k_l}v(\bk):=v(\bk)-v(\bk-\bde_l)$ and $\bar\Delta_\bk^\gamma:=\prod_{l=1}^n\bar\Delta_{k_l}^{\gamma_l}$.} Let $\gamma\in\N_0^n$ with $|\gamma|=n+j$, $j\in\{0,1\}$.  It is easy to see that
   \begin{equation}\label{Eq_Delta_ei_k_eta}
     (e^{-i\eta}-1)^\gamma e^{i\bk\cdot \eta} = (-1)^{|\gamma|} {\bar\Delta}\,_\bk^\gamma e^{i\bk\cdot\eta}.
   \end{equation}
   We will also use the summation by parts formula which states that for $f,h\colon\Z^n\to L(E)$ {of at most polynomial growth, at least one of which is rapidly decreasing,} we have
   \begin{equation}\label{Eq_summ_by_parts_formula}
    \sum_{\bk\in\Z^n} f(\bk) \big[\Delta_{\bk}^\gamma h\big](\bk) = (-1)^{|\gamma|} \sum_{\bk\in\Z^n} \big[ {\bar\Delta}\,_\bk^\gamma f\big](\bk) h(\bk),
   \end{equation}
    see \cite[Proposition~3.1]{ruzhansky-turunen10a}.

  By this, we obtain
  \begin{align*}
    (e^{-i\eta}-1)^\gamma & K_\epsilon(x,\eta,\mu) = (e^{ {-i\eta}}-1)^\gamma \sum_{\bk\in\Z^n} e^{i\bk\cdot\eta} \chi_\epsilon(\bk, \mu) b_\mu(x,x-\eta,\bk)\\
    & = \sum_{\bk\in\Z^n} (-1)^{|\gamma|} \big[ {\bar\Delta}\,_\bk^\gamma e^{i\bk\cdot\eta}\big] \chi_\epsilon(\bk, \mu) b_\mu(x,x-\eta,\bk) \\
    & = \sum_{\bk\in\Z^n} e^{i\bk\cdot\eta} \Delta_\bk^\gamma\big[ \chi_\epsilon(\bk, \mu) b_\mu(x,x-\eta,\bk) \big]\\
    & = \sum_{\bk\in\Z^n} \big( e^{i\bk\cdot\eta} -1\big)      \Delta_\bk^\gamma\big[ \chi_\epsilon(\bk, \mu) b_\mu(x,x-\eta,\bk) \big].
      \end{align*}
    In the last step, we used the equality
    \[ \sum_{\bk\in\Z^n}  \Delta_\bk^\gamma [ \chi_\epsilon(\bk,\mu) b_\mu(x,x-\eta,\bk)  ] = 0\]
    which holds as $\chi_\epsilon(\cdot,\mu)b_\mu(x,x-\eta,\cdot) \in \mathscr S(\Z^n,L(E))$ and all terms appearing in the discrete derivatives cancel.

     Now note that for $\bk\in\Z^n$, $\eta\in\T^n$ and $\theta\in (0,1)$ the elementary inequality
    \[ |e^{i {\bk}\cdot\eta} - 1| \le 2  {|\bk|^\theta } |\eta|^\theta\]
    holds. From this and Lemma~\ref{Lemma 2.6} $b)$, we obtain
    \begin{multline*}
    \big\| (e^{-i\eta}-1)^\gamma   K_\epsilon(x,\eta,\mu) \big\|_{L(E)}  \le \bar{C} |\eta|^{\theta} \bigg[\sum_{\bk\in\Z^n} |\bk|^\theta \langle\bk,\mu\rangle^{-m-\vert\gamma\vert} \\
     + \sum\limits_{\bk\in\Z^n}\vert\bk\vert^\theta\langle \bk,\mu\rangle^{-2m}\langle\bk\rangle^{m-\vert\gamma\vert}\bigg].
    \end{multline*}
     Choosing $\theta:= \theta_j$ for $|\gamma|=n+j,\, j\in\{0,1\}$, as in Lemma~\ref{Lemma2.7}, the last sums can be estimated by $C \mu^{\theta_j}\mu^{-m-j}$.

For completeness, the second sum requires a separate estimate.
For $\mu\ge1$, split it into $0<|\bk|\le\mu$ and $|\bk|>\mu$.
Since $m+\theta_j-j>0$ and $\theta_j<m+j$,
\begin{align*}
\sum_{0<|\bk|\le\mu}|\bk|^{\theta_j}
 \langle\bk,\mu\rangle^{-2m}\langle\bk\rangle^{m-n-j}
&\le C\mu^{-2m}\mu^{m+\theta_j-j},\\
\sum_{|\bk|>\mu}|\bk|^{\theta_j}
 \langle\bk,\mu\rangle^{-2m}\langle\bk\rangle^{m-n-j}
&\le C\mu^{\theta_j-m-j}.
\end{align*}
The range $\mu_0\le\mu\le1$, when nonempty, follows by increasing
the constant.

 By Remark~\ref{Remark 4.7 JFAA}, we see that
    \begin{align*}
     |\eta|^{n+j} \|K_\epsilon(x,\eta,\mu)\|_{L(E)} &\le C\sum_{|\gamma|=n+j} \| (e^{-i\eta}-1)^\gamma K_\epsilon(x,\eta,\mu)\|_{L(E)}\\
      & \le C |\mu\eta|^{\theta_j} \mu^{-m-j}.
      \end{align*}
      Therefore
      \[ \vert\eta\vert^{n}(1+\mu\vert\eta\vert) \|K_\epsilon(x,\eta,\mu)\|_{L(E)} \le C\mu^{-m}\big( \vert\mu\eta\vert^{\theta_0} + \vert\mu\eta\vert^{\theta_1}\big), \]
      that is, \eqref{4-11} holds.

The change of variables $z=\mu\eta$ gives the uniform bound
\[
\int_{\T^n}g_\mu(\eta)\,\dbar\eta
\le C\int_{\R^n}\frac{|z|^{\theta_0}+|z|^{\theta_1}}
 {|z|^n(1+|z|)}\,dz<\infty
\]
since $\theta_0,\theta_1\in (0,1)$. This proves \eqref{4-12}.

    We have seen above that, for fixed $(x,\eta,\mu)$ {and $|\gamma|=n$}, we have
    \[ (e^{-i\eta}-1)^\gamma K_\epsilon(x,\eta,\mu) = \sum_{\bk\in\Z^n} {x_{\bk,\epsilon}(x,\eta,\mu)}\]
    with $x_{\bk,\epsilon}(x,\eta,\mu) := (e^{i\bk\cdot\eta} -1) \Delta_{\bk}^\gamma [ \chi_\epsilon (\bk,\mu)b_\mu(x,x-\eta,\bk)]$.  From the calculations above, for $0<\theta<\min\{1,m\}$, we have the uniform estimate (with respect to $\epsilon$)
    \[ \|x_{ {\bk,\epsilon}}(x,\eta,\mu)\|_{L(E)}\le C y_{\bk}(\eta,\mu),\]
    where
    \[ y_\bk(\eta,\mu) :=   |\eta|^{\theta}  |\bk|^\theta \big[ \langle \bk,\mu\rangle^{-m-n} + \langle\bk,\mu\rangle^{-2m} \langle\bk\rangle^{m-n}\big]\quad (\bk\in\Z^n).\]
     As $\sum_{\bk\in\Z^n} y_\bk(\eta,\mu) <\infty$, the sequence $(y_{\bk}(\eta,\mu))_{\bk\in\Z^n}$ serves as a summable dominating sequence. Because of $\chi_\epsilon(\bk,\mu)\to 1\;(\epsilon\searrow 0)$ for every fixed $\bk$ and $\mu$, we get by dominated convergence the existence of
    \[ \tilde K(x,\eta,\mu) := \lim_{\epsilon\searrow 0} (e^{-i\eta}-1)^\gamma K_\epsilon(x,\eta,\mu)\]
    for every $x$, $\eta$ and $\mu$.

For $\eta\ne0$, choose an index $j$ with $\eta_j\ne0$ and
take $\gamma=n\delta_j$. Then $(e^{-i\eta}-1)^\gamma\ne0$, and division
gives the limit $K(x,\eta,\mu)=\lim_{\epsilon\searrow0}K_\epsilon(x,\eta,\mu)$.
This limit is independent of the chosen index. Set $K(x,0,\mu)=0$.
For fixed $\mu$, every regularized kernel is continuous in $(x,\eta)$, so that its pointwise almost everywhere limit is strongly measurable.

    As the right-hand side of \eqref{4-11} is a dominating and integrable function which is independent of $\epsilon$, we see that $$\|K(x,\cdot,\mu)\|_{L^1(\T^n,L(E))} = \lim_{\epsilon\to 0} \|K_\epsilon(x,\cdot,\mu)\|_{L^1(\T^n,L(E))}$$
    by dominated convergence again, which yields \eqref{4-12} for $K$ instead of $K_\epsilon$.\ {The same argument also gives $K_\epsilon(x,\cdot,\mu)\to K(x,\cdot,\mu)$ in $L^1(\T^n,L(E))$.}
\end{proof}

\begin{theorem}\label{Theorem2.10}
 Let $\mu_0>0$,  $\rho\ge n+1$, $1\le p<\infty$, $k\in\N_0$, $r\geq \max\big\{k,\ \lfloor -m\rfloor +n + {1}\big\}$,  and $\big\{b_\mu\,:\,\mu\ge\mu_0 \big\}\subset S^{-m,\rho,r}(\T^{2n}\times\Z^n,L(E))$  as in Lemma~\ref{Lemma 2.6}\,b). Then $b_\mu(X,Y,D)\in L(W_p^k(\T^n,E))$ and
 \[ \Vert b_\mu(X,Y,D) \Vert_{L(W_p^k(\T^n,E))} \le C\mu^{-m} \qquad (\mu\ge\mu_0) \]
with a constant $C$ independent of $\mu$.
 \end{theorem}

\begin{proof}
 (i) Consider first the case $k=0$. For $f\in\CTE$ we have
 \begin{align*}
  \big(b_\mu(X,Y, D)f\big)(x) & = \ossumint e^{i\bk\cdot y}b_\mu(x,x-y,\bk)f(x-y)\,\dbar y\\
  & = \lim\limits_{\epsilon\searrow 0} \int_{\T^n}\Big( \sum\limits_{\bk\in\Z^n} e^{i\bk\cdot y}\chi_\epsilon(\bk,\mu)b_\mu(x,x-y,\bk)  \Big)f(x-y)\,\dbar y\\
  & = \lim\limits_{\epsilon\searrow 0}\int_{\T^n} K_\epsilon(x,y,\mu)f(x-y)\,\dbar y
 \end{align*}
with $K_\epsilon$ as in Lemma~\ref{Lemma2.9}. By this Lemma, $K_\epsilon(x,\cdot,\mu)\in L^1(\T^n,L(E))$ with $\epsilon$-independent dominating function, and $K_\epsilon$ converges pointwise to a strongly measurable function $K$ for $\epsilon\searrow 0$. By the dominated convergence theorem, we have
\[ \big(b_\mu(X,Y,D)f\big)(x) = \int_{\T^n}K(x,y,\mu)f(x-y)\,\dbar y \qquad (x\in\T^n).\]

{The pointwise estimate $\|K(x,y,\mu)\|\le C\mu^{-m}g_\mu(y)$ and Young\textquotesingle s inequality for convolutions give, for $\mu\ge\mu_0$,}
\begin{align*}
\big\Vert b_\mu(X,Y,D)f \big\Vert_{L^p(\T^n,E)} & \le {C\mu^{-m}\Vert g_\mu\Vert_{L^1(\T^n)}}\Vert f\Vert_{L^p(\T^n,E)}\\
& \le \bar{C}\mu^{-m}\Vert f\Vert_{L^p(\T^n,E)}
\end{align*}
for all $f\in\CTE$.

(ii) Now let $k\in\N$. For $\beta\in\N_0^n$ with $\vert\beta\vert\le k$ and $f\in\CTE$, we first
differentiate the regularized integral in (i), taking the limit,
 and write $D_{x,y}^{\beta'}:=(\partial_x+\partial_y)^{\beta'}$ to obtain
\[
\sum_{\beta'\le\beta}\binom{\beta}{\beta'}
\int_{\T^n}\sum_{\bk\in\Z^n}e^{i\bk\cdot y}\chi_\epsilon(\bk,\mu)
(D_{x,y}^{\beta'}b_\mu)(x,x-y,\bk)
(\partial_x^{\beta-\beta'}f)(x-y)\,\dbar y.
\]
The estimates in Lemma~\ref{Lemma 2.6}\,b) allow
Lemma~\ref{Lemma2.9} to be applied to each differentiated amplitude, and dominated convergence and Young's inequality show that
these expressions converge in $L^p$ as $\epsilon\searrow0$.
The limits are the distributional derivatives of $b_\mu(X,Y,D)f$.\
We thus find
\begin{align*}
 \big\Vert \partial_x^\beta [  b_\mu(X,Y,D)f] \big\Vert_{L^p(\T^n,E)} & \le C\mu^{-m}\sum\limits_{\substack{\gamma\in\N_0^n\\\vert\gamma\vert\le k}}\Vert \partial^\gamma f\Vert_{L^p(\T^n,E)}\\
 & \le C\mu^{-m}\Vert f\Vert_{W_p^k(\T^n,E)}.
\end{align*}
For $1\le p<\infty$, $\CTE$ is dense in $W_p^k(\T^n,E)$ which {gives the unique bounded extension $b_\mu(X,Y,D)\in L(W_p^k(\T^n,E))$ and the asserted norm estimate}.
\end{proof}

\begin{corollary} \label{Corollary_3.11}
 In the situation of Lemma~\ref{Lemma2.5}, {assume additionally $\rho\ge n+1$ and $1\le p<\infty$. Then} we have $a_\lambda^\#(Y,D) \in L(W_p^k(\T^n,E))$ and
 \[ \big\Vert a_\lambda^\#(Y,D)\big\Vert_{L(W_p^k(\T^n,E))} \le C\mu^{-m}\]
 for all $\mu=\vert \lambda\vert^{1/m}\ge\mu_0>0$, {$\Re\lambda\ge0$,} and $k,\,r\in\N_0$ with $r\geq \max\big\{k,\ \lfloor -m\rfloor +n + {1}\big\}$.
\end{corollary}

By Remark~\ref{rem:sector-extension}, the same estimate as in Corollary \ref{Corollary_3.11} holds in
the enlarged sector for sufficiently large $|\lambda|$, with
$a_\lambda^\#:=a_\lambda^{-1}$.

\begin{lemma}\label{Lemma2.12}
 Let $m>0$. For all $0<\theta<1$, $\bar{\theta}>0$ with $\bar{\theta}+\theta < \min\{1,m\}$, we have
 \begin{equation}\label{Eq2-12}
  \sum\limits_{\bk\in\Z^n}\vert \bk\vert^\theta \langle \bk, \mu\rangle^{-m}\langle \bk\rangle^{m-n-1}\le C_\theta\mu^{-\bar{\theta}}
 \end{equation}
for all $\mu\ge\mu_0>0$.
\end{lemma}
\begin{proof}
 {Inequality} \eqref{Eq2-12} is equivalent to
 \[ \sum\limits_{\bk\in\Z^n} \mu^{\bar{\theta}}\vert \bk\vert^\theta \langle \bk, \mu\rangle^{-m}\langle \bk\rangle^{m-n-1}\le C_\theta. \]
{In the following estimates, the zero summand is omitted whenever a negative power of $|\ell|$ occurs, the original summand at $\ell=0$ equals zero.}
Now,
\begin{align*}
 \sum\limits_{\bk\in\Z^n} \mu^{\bar{\theta}} & \vert \bk\vert^\theta \langle \bk, \mu\rangle^{-m}\langle \bk\rangle^{m-n-1}\\
 & = \sum\limits_{\ell\in\mu^{-1}\Z^n} \mu^{\bar{\theta}+\theta}\vert\ell\vert^\theta(1+\mu^2\vert\ell\vert^{2}+\mu^2)^{-\frac{m}{2}}(1+\mu^2\vert\ell\vert^{2})^{\frac{m-n-1}{2}}.
\end{align*}
On the one hand,
\begin{align*}
 \sum\limits_{\substack{\ell\in\mu^{-1}\Z^n\\\vert\ell\vert_\infty >1}} & \mu^{\bar{\theta}+\theta}  \vert\ell\vert^\theta(1+\mu^2\vert\ell\vert^{2}+\mu^2)^{-\frac{m}{2}}(1+\mu^2\vert\ell\vert^{2})^{\frac{m-n-1}{2}}\\
 & \le \sum\limits_{\substack{\ell\in\mu^{-1}\Z^n\\\vert\ell\vert_\infty >1}} \mu^{\bar{\theta}+\theta}\vert\ell\vert^\theta(1+\mu^2\vert\ell\vert^{2})^{-\frac{m}{2}}(1+\mu^2\vert\ell\vert^{2})^{\frac{m-n-1}{2}}\\
 & =  \sum\limits_{\substack{\ell\in\mu^{-1}\Z^n\\\vert\ell\vert_\infty >1}} \mu^{\bar{\theta}+\theta}\vert\ell\vert^\theta(1+\mu^2\vert\ell\vert^{2})^{\frac{-n-1}{2}}   \le  \sum\limits_{\substack{\ell\in\mu^{-1}\Z^n\\\vert\ell\vert_\infty >1}} \mu^{\bar{\theta}+\theta}\vert\ell\vert^\theta \mu^{-n-1}\vert\ell\vert^{-n-1}\\
 & = \sum\limits_{\substack{\ell\in\mu^{-1}\Z^n\\\vert\ell\vert_\infty >1}} \mu^{\bar{\theta}+\theta - 1}\mu^{-n}\vert\ell\vert^{\theta-n-1}  = \mu^{\bar{\theta}+\theta - 1}\sum\limits_{\substack{\ell\in\mu^{-1}\Z^n\\ \vert\ell\vert_\infty >1}}\vert\ell\vert^{\theta-n-1}\mu^{-n}\\
& \!\!\overset{\mbox{\tiny \eqref{2-7}}}{\le} C_\theta\mu^{\bar{\theta}+\theta - 1}  \le C_{\theta,\mu_0}.
\end{align*}
On the other hand, if $m<1$, then $ (1+\mu^2\vert\ell\vert^{2})^{\frac{m-1}{2}} \le 1$ and therefore
 \begin{align*}
  \sum\limits_{\substack{\ell\in\mu^{-1}\Z^n\\\vert\ell\vert_\infty \le 1}} & \mu^{\bar{\theta}+\theta}  \vert\ell\vert^\theta\underbrace{(1+\mu^2\vert\ell\vert^{2}+\mu^2)^{-\frac{m}{2}}}_{\le \mu^{-m}}\;(1+\mu^2\vert\ell\vert^{2})^{\frac{m-n-1}{2}}\\
  & \le \mu^{\bar{\theta}+\theta - m} \sum\limits_{\substack{\ell\in\mu^{-1}\Z^n\\\vert\ell\vert_\infty \le 1}}\vert\ell\vert^\theta(1+\mu^2\vert\ell\vert^{2})^{\frac{m-1}{2}}\underbrace{(1+\mu^2\vert\ell\vert^{2})^{-\frac{n}{2}}}_{\le \mu^{-n}\vert\ell\vert^{-n}}\\
  & \le \mu^{\bar{\theta}+\theta - m}\sum\limits_{\substack{\ell\in\mu^{-1}\Z^n\\\vert\ell\vert_\infty \le 1}}\vert\ell\vert^{\theta-n}\mu^{-n}  \!\!\overset{\mbox{\tiny \eqref{2-6}}}{\le} C_\theta\mu^{\bar{\theta}+\theta - m}  \le C_{\theta,\mu_0}.
 \end{align*}
If $m\ge 1$ and $\vert\ell\vert_\infty \le 1$, then
\[ (1+\mu^2\vert\ell\vert^{2})^{\frac{m-1}{2}} \le C(1+\mu^2\vert\ell\vert_\infty^{2})^{\frac{m-1}{2}} \le C(1+\mu^2)^{\frac{m-1}{2}} = C\langle \mu\rangle^{m-1}.\]
Therefore, if $m\ge 1$, then
\begin{align*}
 \sum\limits_{\substack{\ell\in\mu^{-1}\Z^n\\\vert\ell\vert_\infty \le 1}} & \mu^{\bar{\theta}+\theta}  \vert\ell\vert^\theta\underbrace{(1+\mu^2\vert\ell\vert^{2}+\mu^2)^{-\frac{m}{2}}}_{\le \langle\mu\rangle^{-m}}\;(1+\mu^2\vert\ell\vert^{2})^{\frac{m-n-1}{2}}\\
 & \le \sum\limits_{\substack{\ell\in\mu^{-1}\Z^n\\\vert\ell\vert_\infty \le 1}}\mu^{\bar{\theta}+\theta}\vert\ell\vert^\theta\langle\mu\rangle^{-m}\underbrace{(1+\mu^2\vert\ell\vert^{2})^{-\frac{n}{2}}}_{\le \mu^{-n}\vert\ell\vert^{-n}}\;\underbrace{(1+\mu^2\vert\ell\vert^{2})^{\frac{m-1}{2}}}_{\le C\langle\mu\rangle^{m-1}}\\
 & \le C\sum\limits_{\substack{\ell\in\mu^{-1}\Z^n\\\vert\ell\vert_\infty \le 1}}\mu^{\bar{\theta}+\theta}\vert\ell\vert^{\theta-n}\underbrace{\langle\mu\rangle^{-1}}_{\le \mu^{-1}}\,\mu^{-n}  \le C\mu^{\bar{\theta}+\theta-1}\sum\limits_{\substack{\ell\in\mu^{-1}\Z^n\\\vert\ell\vert_\infty \le 1}}\vert\ell\vert^{\theta-n}\mu^{-n}\\
 & \!\!\overset{\mbox{\tiny \eqref{2-6}}}{\le} CC_\theta\mu^{\bar{\theta}+\theta-1} \le C_{\theta,\mu_0}.
\end{align*}
The estimates above yield the assertion.
\end{proof}

\begin{lemma}\label{lem:remainder-estimate}
Let {$\rho\ge n$}, $\mu_0>0$,
$k\in\mathbb N_0$, $r\geq \max\big\{k,\ \lfloor m\rfloor +n + {1}\big\}$, and $1\le p<\infty$.

\begin{enumerate}

\item[\textnormal{$a)$}]
Let $j\in\{1,\ldots,n\}$ and assume $\rho\ge n+1$.
Let $b_\mu^{(j)}(x,\eta,\bk)$ be continuous for
$(x,\eta)\in\T^n\times[-\pi,\pi)^n$ with continuous derivatives
in $x$ up to order $k$, periodic in $x$, and satisfying
\[
\bigl\|
\Delta_{\bk}^\alpha
{\partial_x^{\beta}}
b_{\mu}^{(j)}(x,\eta,\bk)
\bigr\|_{L(E)}
\le
C_{\alpha,\beta}
\langle \bk,\mu\rangle^{-m}
\langle \bk \rangle^{m-|\alpha|}
\]
for all {$x\in\mathbb T^n$, $\eta\in[-\pi,\pi)^n$}, $\bk\in\mathbb Z^n$, $\mu\ge\mu_0$,
$|\alpha|\le \rho$, and {$|\beta|\le k$} with constants $C_{\alpha, \beta}$ independent of
$x,\eta,\bk$ and $\mu$.\\
Define
\[
(T_{\mu,j}f)(x) := \lim_{\varepsilon\searrow0} \int_{\mathbb T^n} K_{\varepsilon}^{(j)}(x,\eta,\mu) f(x-\eta) \, \dbar\eta,
\]
where
\[
K_{\varepsilon}^{(j)}(x,\eta,\mu) := \sum_{\bk\in\mathbb Z^n} e^{i\bk\cdot\eta} \eta_j \chi_\varepsilon(\bk,\mu) b_{\mu}^{(j)}(x,\eta,\bk)
\]
and $\chi_\varepsilon$ {is} as in Lemma \ref{Lemma 2.6}.

\noindent Then, there exist constants $C>0$ and $\vartheta>0$, independent of $\mu$, such
that
\[
\|T_{\mu,j}\|_{L(W_p^k(\mathbb T^n,E))} \le C\mu^{-\vartheta} \qquad (\mu\ge\mu_0).
\]

\item[\textnormal{$b)$}]
Assume that
\[
\{q_\mu:\mu\ge\mu_0\} \subset S^{m-1,\rho,r} \big(\mathbb T^{2n}\times\mathbb Z^n,L(E)\big)
\]
is bounded and satisfies
\[
\bigl\| \Delta_{\bk}^\alpha \partial_{x,y}^\beta q_\mu(x,y,\bk) \bigr\|_{L(E)} \le C_{\alpha,\beta} \langle \bk,\mu\rangle^{-m} \langle \bk\rangle^{m-1-|\alpha|}
\]
for all $x,y\in\mathbb T^n$, $\bk\in\mathbb Z^n$, $\mu\ge\mu_0$,
$|\alpha|\le \rho$, and $|\beta|\le r$ with constants $C_{\alpha, \beta}$ independent of
$x,y,\bk$ and $\mu$.

\noindent Then, there exist constants $C>0$ and $\vartheta>0$, independent of $\mu$, such
that
\[
\|q_\mu(X,Y,D)\|_{L(W_p^k(\mathbb T^n,E))} \le C\mu^{-\vartheta} \qquad (\mu\ge\mu_0).
\]

\item[\textnormal{$c)$}]
Let $\phi:\mathbb Z^n\to(0,\infty)$ and assume that
\[
\left\{\phi(\bh)^{-1}q_{\mu,\bh}: \mu\ge\mu_0,\ \bh\in\mathbb Z^n \right\}
\]
is bounded in
\[
S^{m-1,\rho,r} \big(\mathbb T^{2n}\times\mathbb Z^n,L(E)\big).
\]
Assume moreover that
\[
\bigl\| \Delta_{\bk}^\alpha \partial_{x,y}^\beta q_{\mu,\bh}(x,y,\bk) \bigr\|_{L(E)}
\le
C_{\alpha,\beta} \phi(\bh) \langle \bk,\mu\rangle^{-m} \langle \bk\rangle^{m-1-|\alpha|}
\]
for all $x,y\in\mathbb T^n$, $\bk,\bh\in\mathbb Z^n$, $\mu\ge\mu_0$,
$|\alpha|\le\rho$ and $|\beta|\le r$, with constants $C_{\alpha, \beta}$ independent of
$x,y,\bk,\bh$ and $\mu$.\\
Then, there exist constants $C>0$ and $\vartheta>0$, independent of $\mu$ and
${\bh}$, such that
\[
\|q_{\mu,\bh}(X,Y,D)\|_{L(W_p^k(\mathbb T^n,E))} \le C\mu^{-\vartheta}\phi(\bh) \qquad (\mu\ge\mu_0,\ \bh\in\mathbb Z^n).
\]

\end{enumerate}
\end{lemma}

\begin{proof}
We first prove $a)$ using similar ideas to the proof of Lemma \ref{Lemma2.9} and that of Theorem \ref{Theorem2.10}. Fix $j\in\{1,\ldots,n\}$ and
$f\in C^\infty(\mathbb T^n,E)$. By the elementary estimate
\[
|\eta_j|
\le
C\,|e^{-i\eta_j}-1|
\qquad \big(\eta_j\in[-\pi,\pi)\big),
\]
we have
\[
\|K_{\varepsilon}^{(j)}(x,\eta,\mu)\|_{L(E)}
\le
C
\bigg\| \sum_{\textbf{k}\in\mathbb Z^n} (e^{-i\eta_j}-1) e^{i\textbf{k}\cdot\eta} \chi_\varepsilon({\bk},\mu)b_{\mu}^{(j)}(x,\eta,\textbf{k})\bigg\|_{L(E)}.
\]
Let $\delta_j$ be the $j$-th canonical unit multi-index. Since
\[
(e^{-i\eta_j}-1)e^{i\textbf{k}\cdot\eta}
=
-{\bar\Delta}\,_\textbf{k}^{\delta_j}e^{i\textbf{k}\cdot\eta}
\]
(compare with \eqref{Eq_Delta_ei_k_eta}), summation by parts on $\mathbb Z^n$ gives
\[
\|K_{\varepsilon}^{(j)}(x,\eta,\mu)\|_{L(E)}
\le C \bigg\| \sum_{\textbf{k}\in\mathbb Z^n} e^{i\textbf{k}\cdot\eta} \Delta_{{\bk}}^{\delta_j}\bigl[\chi_\varepsilon({\bk},\mu)b_{\mu}^{(j)}(x,\eta,\textbf{k})\bigr]\bigg\|_{L(E)}.
\]

By the discrete Leibniz formula, the estimates for the cut-off $\chi_\varepsilon$, and the assumptions on
$b_{\mu}^{(j)}$, we obtain, for $|\alpha|\le \rho-1$ and {$|\beta|\le k$},

\begin{equation*}
\bigl\| \Delta_\textbf{k}^\alpha {\partial_x^\beta} \Delta_\textbf{k}^{\delta_j}[\chi_\varepsilon(\textbf{k},\mu)b_{\mu}^{(j)}(x,\eta,\textbf{k})]\bigr\|_{L(E)}
\le
C_{\alpha,\beta} \langle \textbf{k},\mu\rangle^{-m} \langle \textbf{k}\rangle^{m-1-|\alpha|}.
\end{equation*}
Indeed, the additional discrete difference $\Delta_\textbf{k}^{\delta_j}$ is precisely what lowers the order from $m$ to $m-1$.

Now, let $|\alpha|=n$. Since $\rho\ge n+1$, the previous estimate applies to
$\Delta_\textbf{k}^\alpha\Delta_\textbf{k}^{\delta_j}$. Using \eqref{4-10}
and summation by parts once more, we get, for $0<\theta<1$,
\[
\begin{aligned}
|\eta|^n \|K_{\varepsilon}^{(j)}(x,\eta,\mu)\|
&\le C \sum_{|\alpha|=n} \left\|\sum_{\textbf{k}\in\mathbb Z^n}e^{i\bk\cdot\eta}\Delta_{\bk}^\alpha \Delta_{\bk}^{\delta_j} [\chi_\varepsilon b_{\mu}^{(j)}] \right\| \\
&\le C|\eta|^\theta \sum_{\textbf{k}\in\mathbb Z^n} |\textbf{k}|^\theta \langle \textbf{k},\mu\rangle^{-m} \langle \textbf{k}\rangle^{m-n-1}.
\end{aligned}
\]
By Lemma~\ref{Lemma2.12}, for every $0<\theta<1$ and every
$\vartheta>0$ with $ \vartheta+\theta<\min\{1,m\}$, we have
\[
\sum_{\textbf{k}\in\mathbb Z^n} |\textbf{k}|^\theta \langle \textbf{k},\mu\rangle^{-m} \langle \textbf{k}\rangle^{m-n-1} \le C\mu^{-\vartheta}.
\]
Therefore,
\[
\|K_{\varepsilon}^{(j)}(x,\eta,\mu)\| \le C\mu^{-\vartheta} |\eta|^{\theta-n}.
\]
Since $\theta>0$, the function $\eta\mapsto|\eta|^{\theta-n}$ is integrable on $\mathbb T^n$. Thus,
\[
\|K_{\varepsilon}^{(j)}(x,\cdot,\mu)\|_{L^1(\mathbb T^n,L(E))}
\le
C\mu^{-\vartheta},
\]
uniformly in $x$, $\varepsilon$, and $\mu$.\\

The $n+1$ discrete differences of $\chi_\varepsilon b_\mu^{(j)}$
converge pointwise in $\bk$ to the corresponding differences of
$b_\mu^{(j)}$. The summable bound above, and division by a nonzero
factor $(e^{-i\eta}-1)^\alpha$ away from $\eta=0$, give pointwise
almost everywhere convergence of $K_\varepsilon^{(j)}$.
The factor $\eta_j/(e^{-i\eta_j}-1)$ is bounded, with its value at
zero assigned by continuity. The pointwise bound
$C\mu^{-\vartheta}|\eta|^{\theta-n}$ then gives convergence of the
regularized operators on smooth $f$ in $L^p$ by dominated convergence,
and Young's inequality gives

\[
\|T_{\mu,j}f\|_{L^p(\mathbb T^n,E)} \le C\mu^{-\vartheta}\|f\|_{L^p(\mathbb T^n,E)}.
\]
Furthermore, {for each fixed $\mu$, let $K^{(j)}(\cdot,\cdot,\mu)\colon\T^{2n}\to L(E)$ be the strongly measurable limit, so that} $K_\epsilon^{(j)}(x,\eta,\mu)\to K^{(j)}(x,\eta,\mu)$  as $\epsilon \searrow 0$ pointwise almost everywhere. In order to estimate the Sobolev norm, let $|\beta|\le k$. {Differentiating, first at the regularized level and then passing to the limit in distributions as in Theorem~\ref{Theorem2.10}, the identity}
\[
(T_{\mu,j}f)(x) = \int_{\mathbb T^n} K^{(j)}(x,\eta,\mu)f(x-\eta)\, \dbar\eta
\]
gives a finite sum of terms involving derivatives of
$b_{\mu}^{(j)}$ of order at most $|\beta|\le k$ and derivatives of $f$ of
order at most $k$. Since $r\ge k$, the same kernel estimate applies to each
term. Hence
\[
\|\partial_x^\beta T_{\mu,j}f\|_{L^p(\mathbb T^n,E)} \le C\mu^{-\vartheta}\|f\|_{W_p^k(\mathbb T^n,E)}.
\]
Summing over $|\beta|\le k$ gives
\[
\|T_{\mu,j}f\|_{W_p^k(\mathbb T^n,E)} \le C\mu^{-\vartheta}\|f\|_{W_p^k(\mathbb T^n,E)}.
\]
By density, $T_{\mu,j}$ extends to a bounded operator on
$W_p^k(\mathbb T^n,E)$, with
\[
\|T_{\mu,j}\|_{L(W_p^k(\T^n,E))} \le C\mu^{-\vartheta}.
\]

We now prove $b)$. For $0<\varepsilon<1$, define
\[
K_\varepsilon(x,\eta,\mu) := \sum_{\textbf{k}\in\mathbb Z^n} e^{i\textbf{k}\cdot\eta} \chi_\varepsilon({\bk},\mu) q_\mu(x,x-\eta,\textbf{k}).
\]
Then
\[
(q_\mu(X,Y,D)f)(x) = \lim_{\varepsilon\searrow0} \int_{\mathbb T^n} K_\varepsilon(x,\eta,\mu)f(x-\eta)\, \dbar\eta .
\]
Let $|\alpha|=n$. Using  \eqref{Eq_Delta_ei_k_eta} and summation by parts, we obtain
\[
(e^{-i\eta}-1)^\alpha K_\varepsilon(x,\eta,\mu) = \sum_{\textbf{k}\in\mathbb Z^n} e^{i\textbf{k}\cdot\eta} \Delta_\textbf{k}^\alpha [\chi_\varepsilon(\textbf{k},\mu)q_\mu(x,x-\eta,\textbf{k})].
\]
As before, since the sum of a discrete difference over $\mathbb Z^n$
vanishes after the cut-off is inserted, we may write the last series with
$e^{i\textbf{k}\cdot\eta}-1$ and obtain the factor
$|\eta|^\theta |{\bk}|^\theta$. Hence
\[
|\eta|^n \|K_\varepsilon(x,\eta,\mu)\|_{L(E)} \le C|\eta|^\theta \sum_{\textbf{k}\in\mathbb Z^n} |\textbf{k}|^\theta \langle \textbf{k},\mu\rangle^{-m} \langle \textbf{k}\rangle^{m-n-1}.
\]
Lemma~\ref{Lemma2.12} gives
\[
\|K_\varepsilon(x,\eta,\mu)\|_{L(E)} \le C\mu^{-\vartheta} |\eta|^{\theta-n}.
\]
Therefore,
\[
\|K_\varepsilon(x,\cdot,\mu)\|_{L^1(\mathbb T^n,L(E))} \le C\mu^{-\vartheta}.
\]
The same Sobolev argument as above yields
\[
\|q_\mu(X,Y,D)\|_{L(W_p^k(\T^n,E))} \le C\mu^{-\vartheta}.
\]
This proves $b)$.\\
Finally, $c)$ follows from $b)$ applied to $\widetilde q_{\mu,\bh}:=\phi(\bh)^{-1}q_{\mu,\bh}$.
Indeed, the assumptions imply that
\[
\{\widetilde q_{\mu,\bh}:\mu\ge\mu_0,\ \bh\in\mathbb Z^n\}
\]
satisfies the hypotheses of $b)$ uniformly in $\bh$. Hence
\[
\|\widetilde q_{\mu,\bh}(X,Y,D)\|_{L(W_p^k(\T^n,E))} \le C\mu^{-\vartheta},
\]
with $C$ independent of $\mu$ and ${\bh}$. Since
\[
q_{\mu,\bh}(X,Y,D) = \phi(\bh)\widetilde q_{\mu,\bh}(X,Y,D),
\]
we conclude that
\[
\|q_{\mu,\bh}(X,Y,D)\|_{L(W_p^k(\T^n,E))} \le C\mu^{-\vartheta}\phi(\bh).
\]
The regularized operators in all three parts converge strongly on
$W_p^k(\T^n,E)$ as $\varepsilon\searrow0$. Indeed, the previous argument, using dominated
convergence, also applies to the derivatives of smooth functions.
The uniform $W_p^k(\T^n,E)$ bounds and density of $C^\infty(\T^n,E)$ then give strong
convergence in $W_p^k(\T^n,E)$.
\end{proof}

From now on, we will assume the following.
 
\refstepcounter{assumption}
\noindent\textbf{Assumption~(\theassumption).}\label{assumptionA}
 Let $k\in\mathbb N_0$, $m>0$, $1\le p<\infty$,
$\rho\ge n+1$, and
$r\geq \max\big\{k,\ \lfloor m\rfloor+ n + {1}\big\}$.\\   

Note that the last condition in Assumption~(\ref{assumptionA}) implies that $k\leq r$ and that there exists $N\in\mathbb N$ such that
$$m+n < {N} \leq r.$$

\noindent Recall that, for $f\in C^\infty(\T^n,E)$,
\[
(\op[a]f)(x)= \big(a(X,D)f\big)(x)=\sum_{\bk\in \Z^n} e^{i\bk\cdot x}a(x,\bk)\hat{f}(\bk).
\]
Observe that, since $k\leq r$,
\[
 \op[a]f\in \WTE\quad \text{for all }f\in C^\infty(\T^n,E) .
\]
We will consider the closure $A_k\colon D(A_k)\subset \WTE\to \WTE$ of the linear operator
$$\op[a]_k\colon C^\infty(\T^n,E)\subset \WTE\to \WTE, \quad f\mapsto \op[a]f$$
in $\WTE$. For this, we need the following lemma, which also implies that $A_k$ is densely defined.

\begin{lemma}\label{lem:closability}
The operator $\op[a]_k$ is densely defined and closable.
\end{lemma}

\begin{proof}
{Since $1\le p<\infty$,  $C^\infty(\T^n,E)$ is dense in $W_p^k(\T^n,E)$. We use the bilinear pairing between $E$ and its topological dual $E^\prime$.}
For fixed $\psi\in C^\infty(\mathbb T^n,E')$, define $T_\psi: C^\infty(\mathbb T^n,E)\to\mathbb C$ by

\[
T_\psi(f)
:=
\int_{\mathbb T^n}
\langle \operatorname{op}[a]f(x),\psi(x)\rangle \,\dbar x.
\]

Using the definition of $\operatorname{op}[a]$,

\[
T_\psi(f)
=
\int_{\mathbb T^n}
\left\langle
\sum_{\mathbf{k}\in\mathbb Z^n}
e^{i\mathbf{k}\cdot x}
a(x,\mathbf{k})
\widehat f(\mathbf{k}),
\psi(x)
\right\rangle
\,\dbar x.
\]

Since $f\in C^\infty(\mathbb T^n,E)$, its Fourier coefficients decay rapidly:

\[
\|\widehat f(\mathbf{k})\| \le C_N \langle \mathbf{k}\rangle^{-N} \qquad \forall N \in \N_0.
\]

Therefore, the above series converges absolutely and uniformly{, so that we can} interchange summation and integration:
\[
T_{\psi}(f) = \sum_{\mathbf{k}\in\mathbb Z^n} \int_{\mathbb T^n} e^{i\mathbf{k}\cdot x} \left\langle a(x,\mathbf{k}) \widehat f(\mathbf{k}), \psi(x) \right\rangle \,\dbar x.
\]

Using the transpose operator, $\langle a(x,\mathbf{k})u,\psi(x)\rangle=\langle u,a(x,\mathbf{k})'\psi(x)\rangle$,
we obtain

\[
T_\psi(f) = \sum_{\mathbf{k}\in\mathbb Z^n} \left\langle \widehat f(\mathbf{k}), a_\psi(\mathbf{k}) \right\rangle,
\]

where

\[
a_{\psi}(\mathbf{k}) := \int_{\mathbb T^n} e^{i\mathbf{k}\cdot x} a(x,\mathbf{k})' \psi(x) \,\dbar x \in E'.
\]

Now, let
\[
g_{\bk}(x):=a(x,\bk)'\psi(x).
\]

Choose $N\in\N$ with $m+n<N\le r$, as above. Using \eqref{Eq-L_ky}, integration by parts gives

\begin{align*}
a_\psi(\bk) & = \int_{\mathbb T^n} e^{i\bk\cdot x}
(L_{-\bk,x}^t)^{N}g_{\bk}(x) \,\dbar x.
\end{align*}

Therefore, by \eqref{Eq-L_kytN} it follows that
\[
\|a_\psi(\bk)\|_{E'}
\le
C_r
\langle \bk\rangle^{{-N}}
\sum_{|\alpha|\le r}
\sup_{x\in\mathbb T^n}
\|\partial_x^\alpha g_{\bk}(x)\|_{E'}.
\]

Since $g_{\bk}(x)=a(x,\bk)'\psi(x)$, Leibniz' rule gives
\[
\|\partial_x^\alpha g_{\bk}(x)\|_{E'}
\le
\sum_{\beta\le\alpha}
\binom{\alpha}{\beta}
\|(\partial_x^\beta a(x,\bk))'\|_{L(E')}
\,
\|\partial_x^{\alpha-\beta}\psi(x)\|_{E'}.
\]
Since
\[
\|(\partial_x^\beta a(x,\bk))'\|_{L(E')}
=
\|\partial_x^\beta a(x,\bk)\|_{L(E)},
\]
and by the symbolic estimate
\[
\|\partial_x^\beta a(x,\bk)\|_{L(E)}
\le
C_\beta\langle \bk\rangle^m,
\]
we obtain
\[
\|\partial_x^\alpha g_{\bk}(x)\|_{E'}
\le
C_{\psi,\alpha}\langle \bk\rangle^m.
\]
Consequently,
\[
\|a_\psi(\bk)\|_{E'}
\le
C_{\psi,r}
\langle \bk\rangle^{m-{N}}.
\]

On the other hand, since $f\in W_p^k(\mathbb T^n,E)$,
\[
\|\widehat f(\bk)\|
\le
C\langle \bk\rangle^{-k}
\|f\|_{W_p^k(\mathbb T^n,E)}.
\]
Therefore,
\[
|T_\psi(f)|
\le
\sum_{\bk\in\mathbb Z^n}
\|\widehat f(\bk)\|_E
\|a_\psi(\bk)\|_{E'}\le
C_\psi
\|f\|_{W_p^k(\T^n,E)}
\sum_{\bk\in\mathbb Z^n}
\langle \bk\rangle^{m-{N}-k}.
\]
The series
\[
\sum_{\bk\in\mathbb Z^n}
\langle \bk\rangle^{m-{N}-k}
\]
converges precisely when $m-{N}-k<-n$. Equivalently, ${N}+k>m+n$. Thus, under the condition ${N}>m+n$, there exists a constant $C_\psi>0$ such that
\[
|T_\psi(f)|
\le
C_\psi
\|f\|_{W_p^k(\mathbb T^n,E)}, \hspace{0.3cm} \forall f \in C^\infty(\mathbb T^n,E).
\]
Hence $T_\psi$ extends continuously to $W_p^k(\mathbb T^n,E)$.

Now, let $(f_l)_{l\in \N}\subset C^\infty(\T^n,E)$ with $f_l\to 0$ as $l\to \infty$ in $\WTE$ and $\op[a]f_l\to g$ in $\WTE$ as $l\to \infty$.\ We have to show that $g=0$.
Since $T_\psi$ is continuous, $T_\psi(f_l)\to0$. That is,
\[
\int_{\mathbb T^n}
\langle
\operatorname{op}[a]f_l(x),
\psi(x)
\rangle
\, \dbar x
\to0.
\]
On the other hand, because $\operatorname{op}[a]f_l\to g \quad\text{in }W_p^k(\mathbb T^n,E)$,
and the functional
\[
u\mapsto \int_{\mathbb T^n} \langle u(x),\psi(x)\rangle \,\dbar x
\]
is continuous on $W_p^k(\mathbb T^n,E)$, we also have
\[
\int_{\mathbb T^n}
\langle
\operatorname{op}[a]f_l(x),
\psi(x)
\rangle
\,\dbar x
\to
\int_{\mathbb T^n}
\langle
g(x),
\psi(x)
\rangle
\,\dbar x.
\]
By uniqueness of limits,
\[
\int_{\mathbb T^n}
\langle
g(x),
\psi(x)
\rangle
\,\dbar x
=
0
\]
for every $\psi\in C^\infty(\mathbb T^n,E')$. Hence $g=0$. Therefore, $\operatorname{op}[a]_k$ is closable.
\end{proof}

For $0\le\theta<\pi$ and
$\mu_0>0$, let
\[
\Sigma_{\theta,\mu_0}
:=
\bigl\{
\lambda\in\mathbb C:
|\lambda|\ge {\mu_0^m},\ |\arg\lambda|\le\theta
\bigr\}.
\]
We now state and prove the main result of this paper (Theorem~\ref{Theorem2.13}). For that, we introduce the following additional assumption.\\

\refstepcounter{assumption}
\noindent\textbf{Assumption~(\theassumption).}\label{assumptionB}
There exists a constant $M\in {\N}$ such that
\[
M>n+k+1+|m-1|+\rho
\qquad \text{and} \qquad r\ge k+M.
\]
\begin{theorem}\label{Theorem2.13}
 Under assumptions (\ref{assumptionA}) and (\ref{assumptionB}), let $a\in S^{m,\rho,r}$ be parabolic, i.e., $a\in P_{\omega,\kappa}S^{m,\rho,r}$ for some $\kappa>0$ and $\omega\ge 0$. Then, there exist constants ${C>0}$, $\theta\in (\pi/2,\pi)$ and $\mu_0>0$ such that for the resolvent set $\rho(-A_k)$ of $-A_k$ we have $\Sigma_{\theta,\mu_0}\subset \rho(-A_k)$ and
 \[ \big\Vert (\lambda+A_k)^{-1}\big\Vert_{L(\WTE)} \le {\frac{C}{1+\vert\lambda\vert}} \qquad (\lambda \in \Sigma_{\theta,\mu_0}).\]
 Therefore, $-A_k\colon W_p^k(\T^n,E)\supset D(A_k)\to W_p^k(\T^n,E)$ generates an analytic {$C_0$-semigroup} on $W_p^k(\T^n,E)$.
\end{theorem}

\begin{proof}
Let
\[
a_\lambda(x,\textbf{k}):=\lambda I_E+a(x,\textbf{k}) \qquad \text{and} \qquad \mu:=|\lambda|^{1/m}.
\]
By the parabolicity of $a$ {and the sector extension in Remark~\ref{rem:sector-extension}}, there exist $\theta_0\in(\pi/2,\pi)$ and $\mu_0>0$ such that, for
$\lambda\in\Sigma_{\theta_0,\mu_0}$, the symbol $a_\lambda(x,\textbf{k})$ is
invertible and
\[
a_\lambda^\#(x,\textbf{k}):=a_\lambda(x,\textbf{k})^{-1}
\]
satisfies the estimates of Lemma~\ref{Lemma2.5}{, uniformly in $\arg\lambda$ in this sector. We enlarge $\mu_0$ so that $\mu_0\ge\max\{1,2\omega_0\}$}.

\medskip

\noindent
\textbf{Surjectivity and resolvent estimate.}
Let
\[
B_\lambda^R:=a_\lambda^\#(Y,D).
\]
By Corollary~\ref{Corollary_3.11},
\[
B_\lambda^R\in  L(W_p^k(\mathbb T^n,E))
\]
and
\[
\|B_\lambda^R\|_{ L(W_p^k(\T^n,E))}
\le
C\mu^{-m}.
\]
Since \(\mu^m=|\lambda|\), this is equivalent to
\[
\|B_\lambda^R\|_{ L(W_p^k(\T^n,E))}
\le
\frac{C}{|\lambda|}.
\]
For \(f\in C^\infty(\mathbb T^n,E)\), the composition formula in Lemma \ref{2.3} gives
\[
a_\lambda(X,D)B_\lambda^R f
=
f+P_\lambda f,
\]
where
\[
P_\lambda=P_\lambda(X,Y,D)
\]
has double symbol
\[
P_\lambda(x,y,\textbf{k})
=
\bigl(a_\lambda(x,\textbf{k})-a_\lambda(y,\textbf{k})\bigr)a_\lambda(y,\textbf{k})^{-1}.
\]
Since $a_\lambda=\lambda+a$,
\[
a_\lambda(x,\textbf{k})-a_\lambda(y,\textbf{k})=a(x,\textbf{k})-a(y,\textbf{k}).
\]
Writing $y=x-\eta$, Taylor's formula gives
\[
a(x,\textbf{k})-a(x-\eta,\textbf{k})
=
\sum_{j=1}^n
\eta_j
\int_0^1
\partial_{x_j}a(x-\tau\eta,\textbf{k})\,d\tau .
\]
Thus,
\[
P_\lambda(x,x-\eta,\textbf{k})
=
\sum_{j=1}^n
\eta_j b_{\mu}^{(j)}(x,\eta,\textbf{k}),
\]
where
\[
b_{\mu}^{(j)}(x,\eta,\textbf{k}) := \int_0^1 \partial_{x_j}a(x-\tau\eta,\textbf{k}) a_\lambda(x-\eta,\textbf{k})^{-1} \,d\tau .
\]
Using Lemma~\ref{Lemma2.5} and the symbol estimates for $a${, with $\eta\in[-\pi,\pi)^n$ and derivatives in $x$ up to order $k$,} the family
\(\{b_{\mu}^{(j)}:\mu\ge\mu_0\}\) satisfies the assumptions of
Lemma~\ref{lem:remainder-estimate} $a)${, since $r\ge k+1$}. Hence
\[
\|P_\lambda\|_{ L(W_p^k(\T^n,E))}
\le
C\mu^{-\vartheta}
\]
for some $\vartheta>0$. Taking $\mu_0$ sufficiently large, we obtain
\[
\|P_\lambda\|_{L(W_p^k(\T^n,E))} \le \frac12.
\]
Therefore, $I+P_\lambda$ is invertible on $W_p^k(\mathbb T^n,E)$, and
\[
\|(I+P_\lambda)^{-1}\|_{L(W_p^k(\T^n,E))}
\le
2.
\]
Let $g\in W_p^k(\mathbb T^n,E)$ and define
\[
v:=(I+P_\lambda)^{-1}g, \qquad u:=B_\lambda^R v.
\]
Below, we will verify that $u\in D(A_k)$ and $(\lambda+A_k)u=g$.\ Moreover,
\[
\|u\|_{W_p^k(\T^n,E)} \le \|B_\lambda^R\|_{ L(W_p^k(\T^n,E))} \|(I+P_\lambda)^{-1}\|_{ L(W_p^k(\T^n,E))} \|g\|_{W_p^k(\T^n,E)},
\]
and therefore
\[
\|u\|_{W_p^k(\T^n,E)} \le C\mu^{-m}\|g\|_{W_p^k(\T^n,E)} = \frac{C}{|\lambda|} \|g\|_{W_p^k(\T^n,E)}.
\]
It remains to justify that $u\in D(A_k)$. Choose
\(v_j\in C^\infty(\mathbb T^n,E)\) such that
\[
v_j\to v \quad\text{in }W_p^k(\T^n,E).
\]
Set $u_j:=B_\lambda^R v_j$. Then,
\[
u_j\to u
\quad\text{in }W_p^k(\T^n,E).
\]
The functions $u_j$ are in $D(A_k)$. In fact, let
\[
 d_j(\bk):=\int_{\T^n}e^{-i\bk\cdot y}
 a_\lambda^\#(y,\bk)v_j(y)\,\dbar y.
\]
Integration by parts $r$ times, using the first-order operator in
\eqref{Eq-L_ky}, and Lemma~\ref{Lemma2.5} give
\[
 \|d_j(\bk)\|\le C_{j,\lambda}\langle\bk\rangle^{-m-r}.
\]
Let $u_{j,L}(x):=\sum_{|\bk|\le L}e^{i\bk\cdot x}d_j(\bk)$.
These are smooth trigonometric polynomials. Since $r>n+k$, the Fourier
series for $u_j$ and for $a(X,D)u_j$, differentiated at most $k$ times,
have summable majorants $C_{j,\lambda}\langle\bk\rangle^{k-m-r}$ and
$C_{j,\lambda}\langle\bk\rangle^{k-r}$, respectively. Consequently,
\[
 u_{j,L}\to u_j,\qquad
 \op[a]u_{j,L}\to a(X,D)u_j
 \quad\hbox{in }C^k(\T^n,E).
\]
It follows from Lemma~\ref{2.3} that $a_\lambda(X,D)u_j=(I+P_\lambda)v_j$.
Thus $u_j\in D(A_k)$ and $A_ku_j=a(X,D)u_j$. Therefore,
\[
a_\lambda(X,D)u_j
=
(I+P_\lambda)v_j
\to
(I+P_\lambda)v
=
g
\quad\text{in }W_p^k(\T^n,E).
\]
Since
\[
a_\lambda(X,D)u_j
=
\lambda u_j+a(X,D)u_j,
\]
we obtain
\[
a(X,D)u_j
\to
g-\lambda u
\quad\text{in }W_p^k(\T^n,E).
\]
Because {$A_k$ is closed and $u_j\in D(A_k)$},
it follows that
\[
u\in D(A_k) \quad \text{and} \quad (\lambda+A_k)u=g.
\]
Hence, \(\lambda+A_k\) is surjective.

Furthermore, once injectivity is proved below, the above construction gives
the resolvent estimate
\[
\|(\lambda+A_k)^{-1}\|_{L(W_p^k(\T^n,E))}
\le
C\mu^{-m}
=
\frac{C}{|\lambda|}
\le
{\frac{2C}{1+|\lambda|}}
\]
for \(\lambda\in\Sigma_{\theta_0,\mu_0}\), after enlarging $\mu_0$ if necessary.

\medskip
\noindent
\textbf{Injectivity}.
We now prove that $\ker(\lambda+A_k)=\{0\}$.\ To that end, let
\[
B_\lambda^L:=a_\lambda^\#(X,D).
\]
{Here}, the symbol \(a_\lambda^\#(x,\textbf{k})\) may be regarded as the double symbol
$a_\lambda^\#(x,y,\textbf{k}):=a_\lambda^\#(x,\textbf{k})$, which is independent of $y$.
Thus Theorem~\ref{Theorem2.10}, together with Lemma~\ref{Lemma2.5}, gives
\[
B_\lambda^L\in\ L\big(W_p^k(\mathbb T^n,E)\big).
\]
Let $f\in C^\infty(\mathbb T^n,E)$ and put
\[
g:=a_\lambda(X,D)f.
\]
Since $f$ is smooth and its Fourier coefficients decrease faster than any polynomial, the symbol estimates for $a_\lambda$ imply that
$g\in C^r(\mathbb T^n,E)$. In particular, since $r\ge M$, $g\in C^M(\mathbb T^n,E)$. By Proposition~\ref{Prop. CM-Ck-continuity}, the operator
$B_\lambda^L$ can be applied to $g$, and
\begin{equation}\label{eq:B-lambda-g}
{B_\lambda^L} g(x)
=
\sum_{\bk\in\mathbb Z^n}
e^{i\bk\cdot x}
a_\lambda^\#(x,\bk)\widehat g(\bk),
\end{equation}
where the series and its derivatives of order at most $k$ converge
absolutely and uniformly. In fact, since $M>k$ and convergence in $C^M(\T^n,E)$ implies convergence in $W^k_p(\T^n, E)$, we can approximate $g$ by smooth functions in the topology of $C^M(\T^n,E)$ and see that both extensions for $B_\lambda^L$, given by Proposition~\ref{Prop. CM-Ck-continuity} and Theorem~\ref{Theorem2.10}, converge to the same limit in $W^k_p(\T^n,E)$.

Let $\chi_\varepsilon$ {be} as in Lemma \ref{Lemma 2.6} and define for $\varepsilon>0$:
\begin{equation}\label{eq:I-lambda-epsilon}
I_{\lambda,\varepsilon}f(x)
:=
\sum_{\bk\in\mathbb Z^n}
e^{i\bk\cdot x}
\chi_{\varepsilon}(\bk, \mu)
a_\lambda^\#(x,\bk)\widehat g(\bk).
\end{equation}
By Proposition~\ref{Prop. CM-Ck-continuity} and dominated convergence,
\begin{equation}\label{eq:I-epsilon-limit}
I_{\lambda,\varepsilon}f
\to
{B_\lambda^L} a_\lambda(X,D)f
\end{equation}
in $C^k(\mathbb T^n,E)$, and hence also in
$W_p^k(\mathbb T^n,E)$, as $\varepsilon\searrow0$.

Now,
\[
g(y)
=
\sum_{\boldsymbol{\ell}\in\mathbb Z^n}
e^{i\boldsymbol{\ell}\cdot y}
a_\lambda(y,\boldsymbol{\ell})
\widehat f(\boldsymbol{\ell}),
\]
and therefore
\begin{align}
\widehat g(\bk)
={}&
\sum_{\boldsymbol{\ell}\in\mathbb Z^n}
\int_{\mathbb T^n}
e^{i(\boldsymbol{\ell}-\bk)\cdot y}
a_\lambda(y,\boldsymbol{\ell})
\widehat f(\boldsymbol{\ell})
\,\dbar y.
\label{eq:Fourier-g}
\end{align}
Substitution into \eqref{eq:I-lambda-epsilon} yields
\begin{align}
I_{\lambda,\varepsilon}f(x) = \sum_{\bk,\boldsymbol{\ell}\in\mathbb Z^n}
e^{i\bk\cdot x} {\chi_\varepsilon(\bk,\mu)}a_\lambda^\#(x,\bk) \int_{\mathbb T^n}
e^{i(\boldsymbol{\ell}-\bk)\cdot y}
a_\lambda(y,\boldsymbol{\ell}) \,\dbar y \widehat f(\boldsymbol{\ell}). \label{eq:double-composition}
\end{align}
Since
\begin{align}
a_\lambda(y,\boldsymbol{\ell})
={}&
a_\lambda(x,\bk)
+
\bigl(a(y,\boldsymbol{\ell})-a(y,\bk)\bigr)
\nonumber\\
&+
\bigl(a(y,\bk)-a(x,\bk)\bigr)
\label{eq:exact-symbol-decomposition}
\end{align}
and $a_\lambda^\#(x,\bk)a_\lambda(x,\bk)=I_E$, the first term gives
\[
P_\varepsilon f(x)
:=
\sum_{\bk\in\mathbb Z^n}
e^{i\bk\cdot x}
{\chi_\varepsilon(\bk,\mu)}\widehat f(\bk).
\]
Thus,
\begin{equation}\label{eq:I-decomposition}
I_{\lambda,\varepsilon}f = P_\varepsilon f + Q_{\lambda,\varepsilon}^{(d)}f + Q_{\lambda,\varepsilon}^{(s)}f,
\end{equation}
where $Q_{\lambda,\varepsilon}^{(d)}$ and
$Q_{\lambda,\varepsilon}^{(s)}$ denote, respectively, the discrete and
spatial difference operators.

We first consider the spatial difference operator
\begin{align}
Q_{\lambda,\varepsilon}^{(s)}f(x)
{\coloneqq}&
\sum_{\bk,\boldsymbol{\ell}\in\mathbb Z^n}
e^{i\bk\cdot x}
{\chi_\varepsilon(\bk,\mu)}
a_\lambda^\#(x,\bk)
\nonumber\\
&\quad\times
\int_{\mathbb T^n}
e^{i(\boldsymbol{\ell}-\bk)\cdot y}
\bigl(a(y,\bk)-a(x,\bk)\bigr)
\,\dbar y\,\widehat f(\boldsymbol{\ell}). \label{eq:Qx-two-sums}
\end{align}
Since the factor $a(y,\bk)-a(x,\bk)$ does not depend on $\boldsymbol{\ell}$, Fourier inversion gives
\[
\sum_{\boldsymbol{\ell}\in\mathbb Z^n}
e^{i\boldsymbol{\ell}\cdot y}
\widehat f(\boldsymbol{\ell})
=
f(y).
\]
Consequently, \eqref{eq:Qx-two-sums} becomes
\begin{align}
Q_{\lambda,\varepsilon}^{(s)}f(x)
={}&
\sum_{\bk\in\mathbb Z^n}
e^{i\bk\cdot x}
{\chi_\varepsilon(\bk,\mu)}
\nonumber\\
&\quad\times
\int_{\mathbb T^n}
e^{-i\bk\cdot y}
a_\lambda^\#(x,\bk)
\bigl(a(y,\bk)-a(x,\bk)\bigr)
f(y)\,\dbar y.
\label{eq:Qx-one-sum}
\end{align}
Let $\eta:=x-y$. By the fundamental theorem of calculus,
\begin{align}
a(x-\eta,\bk)-a(x,\bk)
&=
{\int_0^1\frac{d}{dt}a(x-t\eta,\bk)\,dt}
\nonumber\\
&=
-\sum_{j=1}^n
\eta_j
\int_0^1
\partial_{x_j}a(x-t\eta,\bk)\,dt.
\label{eq:spatial-fundamental-theorem}
\end{align}
We define
\begin{equation}\label{eq:b-mu-j}
b_{\mu}^{(j)}(x,\eta,\bk)
:=
-
a_\lambda^\#(x,\bk)
\int_0^1
\partial_{x_j}a(x-t\eta,\bk)\,dt.
\end{equation}
Then,
\begin{equation}\label{eq:spatial-factorization}
a_\lambda^\#(x,\bk)
\bigl(a(x-\eta,\bk)-a(x,\bk)\bigr)
=
\sum_{j=1}^n
\eta_j b_{\mu}^{(j)}(x,\eta,\bk).
\end{equation}
\noindent
{For $|\alpha|\le\rho$ and $|\beta|\le k$, the discrete and continuous
Leibniz rules, Lemma~\ref{Lemma2.5}, and $r\ge k+1$ give}
\begin{equation}\label{eq:b-mu-estimate}
{\left\|\Delta_\bk^\alpha\partial_x^\beta
b_\mu^{(j)}(x,\eta,\bk)\right\|_{L(E)}
\le C_{\alpha,\beta}\langle\bk,\mu\rangle^{-m}
\langle\bk\rangle^{m-|\alpha|},}
\end{equation}
where $x\in\T^n$ and $\eta\in[-\pi,\pi)^n$. 

Therefore, Lemma~\ref{lem:remainder-estimate} $a)$ applies
to each amplitude $b_\mu^{(j)}$. It follows that there exist
$\vartheta_s>0$ and $C>0$, independent of $\lambda$, such that the limit
\[
Q_\lambda^{(s)}f := \lim_{\varepsilon\searrow0}Q_{\lambda,\varepsilon}^{(s)}f
\]
exists in $W_p^k(\mathbb T^n,E)$ and
\begin{equation}\label{eq:Qx-estimate}
\left\| Q_\lambda^{(s)} \right\|_{ L(W_p^k(\T^n,E))} \le C\mu^{-\vartheta_s}.
\end{equation}
We now consider the discrete difference operator:
\begin{align}
Q_{\lambda,\varepsilon}^{(d)}f(x)
{\coloneqq}&
\sum_{\bk,\boldsymbol{\ell}\in\mathbb Z^n}
e^{i\bk\cdot x}
{\chi_\varepsilon(\bk,\mu)}
a_\lambda^\#(x,\bk)
\nonumber\\
&\quad\times
\int_{\mathbb T^n}
e^{i(\boldsymbol{\ell}-\bk)\cdot y}
\bigl(a(y,\boldsymbol{\ell})-a(y,\bk)\bigr)
\,\dbar y\,
\widehat f(\boldsymbol{\ell}).
\label{eq:Qd-two-sums}
\end{align}
Setting $\mathbf h:=\boldsymbol{\ell}-\bk$, we obtain
\begin{align}
Q_{\lambda,\varepsilon}^{(d)}f(x)
={}&
\sum_{\mathbf h\in\mathbb Z^n}
\sum_{\bk\in\mathbb Z^n}
e^{i\bk\cdot x}
{\chi_\varepsilon(\bk,\mu)}
c_{\mu,\mathbf h}(x,\bk)
\widehat f(\bk+\mathbf h),
\label{eq:Qd-h-decomposition}
\end{align}
where
\begin{equation}\label{eq:c-mu-h}
c_{\mu,\mathbf h}(x,\bk) := a_\lambda^\#(x,\bk)\int_{\mathbb T^n}e^{i\mathbf h\cdot y}
\bigl(a(y,\bk+\mathbf h)-a(y,\bk)\bigr)\,\dbar y.
\end{equation}
For each $\mathbf h\in\mathbb Z^n$, define the operator
\[
M_{-\mathbf h}f(x) := e^{-i\mathbf h\cdot x}f(x).
\]
Its Fourier coefficients satisfy $\widehat{M_{-\mathbf h}f}(\bk)=\widehat f(\bk+\mathbf h)$.
Hence, if
\[
T_{\mu,\mathbf h,\varepsilon}
:=
\operatorname{op}
\bigl[
{\chi_\varepsilon(\cdot,\mu)}c_{\mu,\mathbf h}
\bigr],
\]
then
\begin{equation}\label{eq:Qd-single-sum-family}
Q_{\lambda,\varepsilon}^{(d)}f = \sum_{\mathbf h\in\mathbb Z^n}T_{\mu,\mathbf h,\varepsilon}M_{-\mathbf h}f.
\end{equation}
We next estimate $c_{\mu,\mathbf h}$. Let
\[
D_{\mathbf h}a(y,\bk)
:=
a(y,\bk+\mathbf h)-a(y,\bk).
\]
Using a telescopic sum, the estimates for $a$, and Peetre\textquotesingle s inequality, for $|\alpha|\le\widetilde\rho:=\rho-1$ and $|\beta|\le r$,
\begin{equation}\label{eq:discrete-first-difference}
\left\| \Delta_\bk^\alpha \partial_y^\beta D_{\mathbf h}a(y,\bk) \right\|_{L(E)}
\le C_{\alpha,\beta} \phi(\mathbf h) \langle\bk\rangle^{m-1-|\alpha|},
\end{equation}
where $\phi(\mathbf h):=\langle\mathbf h\rangle^{1+|m-1|+\rho}$.

After $M$ integrations by parts on the torus, we obtain {as in} the proof of Proposition \ref{Prop. CM-Ck-continuity} that
\begin{equation} \label{eq:c-mu-h-after-IBP}
c_{\mu,\mathbf h}(x,\bk)=\sum_{|\gamma|\le M}
{c_{\gamma,M}(-\mathbf h)} a_\lambda^\#(x,\bk) \int_{\mathbb T^n}
e^{i\mathbf h\cdot y}\partial_y^\gamma D_{\mathbf h}a(y,\bk) \,\dbar y
\end{equation}
with
\begin{equation}\label{eq:Lh-coefficient-estimate}
\big|{c_{\gamma,M}(-\mathbf h)}
\big|
\le
C_M\langle\mathbf h\rangle^{-M}.
\end{equation}
For $|\gamma|\le M$, we define
\begin{equation}\label{eq:c-mu-h-gamma}
c_{\mu,\mathbf h,\gamma}(x,\bk)
:=
a_\lambda^\#(x,\bk) \int_{\mathbb T^n}e^{i\mathbf h\cdot y}\partial_y^\gamma
D_{\mathbf h}a(y,\bk)\,\dbar y.
\end{equation}
Each $c_{\mu,\mathbf h,\gamma}$ may be regarded as a symbol independent of the second spatial variable.

The discrete Leibniz rule, Lemma \ref{Lemma2.5}, and
\eqref{eq:discrete-first-difference} imply
\begin{equation} \label{eq:c-mu-h-gamma-estimate}
\left\| \Delta_\bk^\alpha \partial_x^\beta c_{\mu,\mathbf h,\gamma}(x,\bk) \right\|_{L(E)}
\le
C_{\alpha,\beta,\gamma} \phi(\mathbf h) \langle\bk,\mu\rangle^{-m} {\langle\bk\rangle^{m-1-|\alpha|}}
\end{equation}
for $|\alpha|\le\widetilde\rho$ and {$|\beta|\le r$}.
Indeed, derivatives in $x$ act only on $a_\lambda^\#(x,\bk)$,
whereas the $|\gamma|\le M\le r$ derivatives in $y$ act only on
$D_{\mathbf h}a(y,\bk)$. 

Now, Lemma~\ref{lem:remainder-estimate} $c)$ applied to the
family $\bigl\{c_{\mu,\mathbf h,\gamma}:\mu\ge\mu_0,\ \mathbf h\in\mathbb Z^n \bigr\}$,
gives constants $\vartheta_d>0$ and $C>0$, independent of
$\mu$, $\mathbf h$, and $\gamma$, such that
\begin{equation}\label{eq:c-h-operator-estimate}
\left\|
\operatorname{op}
[c_{\mu,\mathbf h,\gamma}]
\right\|_{L(W_p^k(\T^n,E))}
\le
C\mu^{-\vartheta_d}\phi(\mathbf h).
\end{equation}
Using \eqref{eq:c-mu-h-after-IBP},
\eqref{eq:Lh-coefficient-estimate}, and
\eqref{eq:c-h-operator-estimate}, we obtain
\begin{equation}\label{eq:T-h-estimate}
\left\|
T_{\mu,\mathbf h}
\right\|_{L(W_p^k(\T^n,E))}
\le
C
\mu^{-\vartheta_d}
\langle\mathbf h\rangle^{-M}
\phi(\mathbf h),
\end{equation}
where $T_{\mu,\mathbf h}:=\operatorname{op}[c_{\mu,\mathbf h}]$ is the oscillatory limit of
$T_{\mu,\mathbf h,\varepsilon}$.

Now, the operator $M_{-\mathbf h}$ satisfies
\begin{equation}\label{eq:modulation-estimate}
\left\|
M_{-\mathbf h}f
\right\|_{W_p^k(\T^n,E)}
\le
C_k
\langle\mathbf h\rangle^k
\|f\|_{W_p^k(\T^n,E)}.
\end{equation}
Indeed, for every $|\alpha|\le k$,
\[
\partial_x^\alpha
\bigl(e^{-i\mathbf h\cdot x}f(x)\bigr)
=
\sum_{\beta\le\alpha}
\binom{\alpha}{\beta}
(-i\mathbf h)^{\alpha-\beta}
e^{-i\mathbf h\cdot x}
\partial_x^\beta f(x).
\]
It follows from \eqref{eq:T-h-estimate} and
\eqref{eq:modulation-estimate} that
\begin{align}
\left\|
T_{\mu,\mathbf h}M_{-\mathbf h}
\right\|_{L(W_p^k(\T^n,E))}
&\le
C
\mu^{-\vartheta_d}
\langle\mathbf h\rangle^{-M+k}
\phi(\mathbf h)
\nonumber\\
&=
C
\mu^{-\vartheta_d}
\langle\mathbf h\rangle^{
-M+k+1+|m-1|+\rho
}.
\label{eq:ThMh-estimate}
\end{align}
By assumption (\ref{assumptionB}), $-M+k+1+|m-1|+\rho<-n$, and therefore
\begin{equation}\label{eq:h-series-converges}
\sum_{\mathbf h\in\mathbb Z^n}
\langle\mathbf h\rangle^{
-M+k+1+|m-1|+\rho
}
<\infty.
\end{equation}
Consequently, the operator series
\begin{equation}\label{eq:Qd-operator-series}
Q_\lambda^{(d)}
:=
\sum_{\mathbf h\in\mathbb Z^n}
T_{\mu,\mathbf h}M_{-\mathbf h}
\end{equation}
converges absolutely in $L\bigl(W_p^k(\mathbb T^n,E)\bigr)$, and
\begin{equation}\label{eq:Qd-final-estimate}
\left\| Q_\lambda^{(d)} \right\|_{L(W_p^k(\T^n,E))} \le C\mu^{-\vartheta_d}.
\end{equation}
For each fixed $\mathbf h$,
Lemma~\ref{lem:remainder-estimate} gives
$T_{\mu,\mathbf h,\varepsilon}v\to T_{\mu,\mathbf h}v$ in $W_p^k(\T^n,E)$
for every $v\in W_p^k(\T^n,E)$  with $T_{\mu,\mathbf h,\varepsilon}$ also satisfying  \eqref{eq:T-h-estimate} and \eqref{eq:ThMh-estimate}. Therefore, by dominated convergence, for each $f\in W_p^k(\T^n,E)$, we get
\begin{equation}\label{eq:Qd-epsilon-limit}
 Q_{\lambda,\varepsilon}^{(d)}f\to Q_\lambda^{(d)}f
 \quad\hbox{in }W_p^k(\T^n,E).
\end{equation}
Since $P_\varepsilon f\to f$ in $W_p^k(\mathbb T^n,E)$, passing to the limit in
\eqref{eq:I-decomposition} gives
\begin{equation}\label{eq:left-parametrix-smooth}
{B_\lambda^L} a_\lambda(X,D)f
= f+ Q_\lambda^{(d)}f+ Q_\lambda^{(s)}f, \qquad f\in C^\infty(\mathbb T^n,E).
\end{equation}
Set
\[
Q_\lambda
:=
Q_\lambda^{(d)}+{Q_\lambda^{(s)}}
\]
and $\vartheta := \min\{\vartheta_d,\vartheta_s\}$.
From \eqref{eq:Qx-estimate} and \eqref{eq:Qd-final-estimate},
\begin{equation}\label{eq:Q-lambda-small}
\left\| Q_\lambda \right\|_{L(W_p^k(\T^n,E))} \le C\mu^{-\vartheta}.
\end{equation}
Taking \(\mu_0\) sufficiently large, we have for $\mu \geq \mu_0$ that
\[
\|Q_\lambda\|_{L(W_p^k(\T^n,E))}<\dfrac{1}{2} \qquad \quad \text{for all} \, \lambda\in\Sigma_{\theta_0,\mu_0}.
\]
Thus $I+Q_\lambda$ is invertible on \(W_p^k(\mathbb T^n,E)\) for all $\lambda\in\Sigma_{\theta_0,\mu_0}$.

Note that Equation \eqref{eq:left-parametrix-smooth} can be written as
\begin{equation}\label{eq:left-parametrix-IQ}
{B_\lambda^L} a_\lambda(X,D)f
=
(I+Q_\lambda)f,
\qquad
f\in C^\infty(\mathbb T^n,E).
\end{equation}
We recall that $A_k$ is the closure in $W_p^k(\mathbb T^n,E)$ of $a(X,D)$ with $D(a(X,D))=C^\infty(\mathbb T^n,E)$.

Let $u\in D(A_k)$. Then there exists a sequence $(f_j)_{j\in\mathbb N}
\subset C^\infty(\mathbb T^n,E)$ such that $f_j\to u$ and
\[
a(X,D)f_j\to A_ku \quad \text{in} \, W_p^k(\mathbb T^n,E).
\]
Therefore,
\[
a_\lambda(X,D)f_j = \lambda f_j+a(X,D)f_j \to (\lambda+A_k)u \quad \text{in} \, W_p^k(\mathbb T^n,E).
\]
Applying \eqref{eq:left-parametrix-IQ} to $f_j$, we obtain
\[
{B_\lambda^L} a_\lambda(X,D)f_j
=
(I+Q_\lambda)f_j.
\]
Since ${B_\lambda^L},Q_\lambda \in  L\bigl(W_p^k(\mathbb T^n,E)\bigr)$, we conclude that
\begin{equation}\label{eq:parametrix-domain}
{B_\lambda^L}(\lambda+A_k)u = (I+Q_\lambda)u, \qquad u\in D(A_k).
\end{equation}
Suppose now that $u\in D(A_k)$ satisfies
\[
(\lambda+A_k)u=0.
\]
By \eqref{eq:parametrix-domain},
\[
u= (I+Q_\lambda)^{-1}(I+Q_\lambda)u= (I+Q_\lambda)^{-1}[{B_\lambda^L}(\lambda+A_k)u]=0.
\]
Therefore,
\[
\ker(\lambda+A_k)=\{0\}, \qquad \lambda\in\Sigma_{\theta_0,\mu_0}.
\]
Thus, $\lambda+A_k$ is injective.

Consequently, \(\lambda+A_k\) is bijective for all
\(\lambda\in\Sigma_{\theta_0,\mu_0}\) with $\mu_0$ sufficiently large. The
resolvent estimate obtained in the surjectivity part now gives
\[
\|(\lambda+A_k)^{-1}\|_{ L(W_p^k(\T^n,E))}
\le
\frac{C}{|\lambda|}
\le
{\frac{2C}{1+|\lambda|}},
\qquad
\lambda\in\Sigma_{\theta_0,\mu_0}.
\]
From this and the dense domain of $A_k$ it follows that $-A_k$ generates an
analytic $C_0$-semigroup on $W_p^k(\mathbb T^n,E)$.
\end{proof}

Now, we can consider the generation of semigroups on Besov spaces.\ Assume the hypotheses of Theorem~\ref{Theorem2.13} for an integer $k\ge1$,
and set $X_0=L^p(\T^n,E)$, $X_k=W_p^k(\T^n,E)$.
For $0<s<k$ and $1\le q\le\infty$, the periodic Besov space is
\begin{equation}\label{eq:Besov-interpolation}
 B_{pq}^s(\T^n,E)=(X_0,X_k)_{s/k,q}
\end{equation}
with equivalent norms. Here $(\cdot,\cdot)_{\theta,q}$ denotes real
interpolation by the $K$-method. Equivalently, for a smooth inhomogeneous
dyadic partition $(\varphi_j)_{j\ge0}$ on $\R^n$, its norm is
\[
 \|f\|_{B_{pq}^s(\T^n,E)}
 =\big\|\big(2^{js}\|\op[\varphi_j]f\|_{L^p(\T^n,E)}\big)_{j\ge0}\big\|_{\ell^q}.
\]
No geometric hypothesis on $E$ is needed in \eqref{eq:Besov-interpolation}.
Indeed,  in the same way as in \cite[Section 3.3]{amann09}, where $\R^n$ is considered 
instead of $\T^n$,  one sees that $$B_{p1}^\ell(\T^n,E)\hookrightarrow W_p^\ell(\T^n,E)\hookrightarrow B_{p\infty}^\ell(\T^n,E)$$
for $\ell=0,k$. Applying real interpolation to  these inclusions  and using the interpolation properties of Besov spaces,  we obtain \eqref{eq:Besov-interpolation}; see
\cite[Section~5]{amann97} for the interpolation of vector-valued Besov spaces 
and \cite{BDHN2018} for the toroidal dyadic setting.

\begin{corollary}\label{cor:Besov-generation}
Under the assumptions of Theorem~\ref{Theorem2.13}, define the part of $A_0$ in $B_{pq}^s(\T^n,E)$ by
\[
 \begin{split}
 D(A_{s,q})&:=\{u\in D(A_0)\cap B_{pq}^s(\T^n,E):
                    A_0u\in B_{pq}^s(\T^n,E)\},\quad A_{s,q}u:=A_0u.
 \end{split}
\]
There exist $\theta>\pi/2$, $\mu_0>0$, and $C>0$ such that
\[
 \Sigma_{\theta,\mu_0}\subset\rho(-A_{s,q}),\qquad
 \|(\lambda+A_{s,q})^{-1}\|_{L(B_{pq}^s(\T^n,E))}
 \le \frac{C}{1+|\lambda|}.
\]
If $q<\infty$, $-A_{s,q}$ generates an analytic $C_0$-semigroup.
For $q=\infty$, the corresponding analytic semigroup on $B_{p\infty}^s(\T^n,E)$
need not be strongly continuous at zero. Its restriction to
\[
 b_{p\infty}^s(\T^n,E):=\overline{C^\infty(\T^n,E)}^{\,B_{p\infty}^s(\T^n,E)}
\]
is an analytic $C_0$-semigroup, whose generator is the part of $-A_0$
in this closed subspace.
\end{corollary}
\begin{proof}
The hypotheses for $k$ imply those for $k=0$. As $A_k$ is defined as the closure of $\op[a]_k$, we obtain $D(A_k)\subset D(A_0)$ and $A_k=A_0$ on $D(A_k)$.
Uniqueness of the resolvent therefore gives
\[
 (\lambda+A_0)^{-1}|_{X_k}=(\lambda+A_k)^{-1}.
\]
The two estimates in Theorem~\ref{Theorem2.13}, on a common sector,
can consequently be interpolated. For $\sigma=s/k$ this yields
\[
 \|(\lambda+A_0)^{-1}\|_{L(B_{pq}^s(\T^n,E))}
 \le \frac{C_{\mathrm{int}}C_0^{1-\sigma}C_k^\sigma}{1+|\lambda|}.
\]
The constant $C_{\mathrm{int}}$ accounts for the equivalent Besov norm.
The restriction maps $B_{pq}^s(\T^n,E)$ into $D(A_{s,q})$, because
$A_0(\lambda+A_0)^{-1}f=f-\lambda(\lambda+A_0)^{-1}f$.
It is thus the resolvent of the stated part.

The analytic semigroups on $X_0$ and $X_k$ are consistent, as follows
from the common resolvent representation. Interpolation gives an analytic
semigroup, locally bounded at zero, on $(X_0,X_k)_{\sigma,q}$.
For $q<\infty$, $X_k$ is dense in this interpolation space, so strong
continuity follows from strong continuity on $X_k$ and uniform boundedness
near zero. For $q=\infty$, the same argument applies to the closure of
$X_k$, which equals $b_{p\infty}^s(\T^n,E)$ by periodic smoothing and the continuous
embedding $X_k\hookrightarrow B_{p\infty}^s(\T^n,E)$. In all cases, the generator of the $C_0$-semigroup is identified by the resolvent and therefore 
given as the part of $-A_0$ in $B^s_{pq}(\T^n,E)$ for $q<\infty$ and in $b^s_{p\infty}(\T^n,E)$ for $q=\infty$.
\end{proof}

\bibliographystyle{abbrv}

\begin{thebibliography}{10}

\bibitem{Abels-Pfeuffer18}
H.~Abels and C.~Pfeuffer.
\newblock Characterization of non-smooth pseudodifferential operators.
\newblock {\em J. Fourier Anal. Appl.}, 24(2):371--415, 2018.

\bibitem{agranovich79}
M.~S. Agranovi{\v c}.
\newblock Spectral properties of elliptic pseudodifferential operators on a
  closed curve.
\newblock {\em Funktsional. Anal. i Prilozhen.}, 13(4):54--56, 1979.

\bibitem{amann97}
H.~Amann.
\newblock Operator-valued {F}ourier multipliers, vector-valued {B}esov spaces,
  and applications.
\newblock {\em Math. Nachr.}, 186:5--56, 1997.

\bibitem{amann09}
H.~Amann.
\newblock {\em Anisotropic function spaces and maximal regularity for parabolic
  problems. {P}art 1}, volume~6 of {\em Jind{\v r}ich Ne{\v c}as Center for
  Mathematical Modeling Lecture Notes}.
\newblock Matfyzpress, Prague, 2009.
\newblock Function spaces.

\bibitem{amann-walker05}
H.~Amann and C.~Walker.
\newblock Local and global strong solutions to continuous
  coagulation-fragmentation equations with diffusion.
\newblock {\em J. Differential Equations}, 218(1):159--186, 2005.

\bibitem{amosov88}
B.~A. Amosov.
\newblock On the theory of pseudodifferential operators on the circle.
\newblock {\em Uspekhi Mat. Nauk}, 43(3(261)):169--170, 1988.

\bibitem{arendt-bu02}
W.~Arendt and S.~Bu.
\newblock The operator-valued {M}arcinkiewicz multiplier theorem and maximal
  regularity.
\newblock {\em Math. Z.}, 240(2):311--343, 2002.

\bibitem{arendt-bu04}
W.~Arendt and S.~Bu.
\newblock Operator-valued {F}ourier multipliers on periodic {B}esov spaces and
  applications.
\newblock {\em Proc. Edinb. Math. Soc. (2)}, 47(1):15--33, 2004.

\bibitem{bahouri-chemin-danchin11}
H.~Bahouri, J.-Y. Chemin, and R.~Danchin.
\newblock {\em Fourier {A}nalysis and {N}onlinear {P}artial {D}ifferential
  {E}quations}, volume 343 of {\em Grundlehren der mathematischen
  {W}issenschaften}.
\newblock Springer, Heidelberg, 2011.

\bibitem{barraza-denk-hernandez-nau}
B.~Barraza~Mart\'inez, R.~Denk, J.~Hern\'andez~Monz\'on, and T.~Nau.
\newblock Generation of semigroups for vector-valued pseudodifferential
  operators on the torus.
\newblock {\em J. Fourier Anal. Appl.}, 22(4):823--853, 2016.

\bibitem{BDHN2018}
B.~Barraza~Mart\'inez, R.~Denk, J.~Hern\'andez~Monz\'on, and M.~Nendel.
\newblock Mapping properties for operator-valued pseudodifferential operators
  on toroidal {B}esov spaces.
\newblock {\em J. Pseudo-Differ. Oper. Appl.}, 9(3):523--538, 2018.

\bibitem{Benhamoud2024Endpoint}
R.~Benhamoud.
\newblock Endpoint boundedness of toroidal pseudo-differential operators.
\newblock {\em Open Mathematics}, 22(1):20240023, 2024.

\bibitem{bu-kim04}
S.~Bu and J.-M. Kim.
\newblock Operator-valued {F}ourier multiplier theorems on {$L_p$}-spaces on
  {$\mathbb{T}^d$}.
\newblock {\em Arch. Math. (Basel)}, 82(5):404--414, 2004.

\bibitem{cardona14}
D.~Cardona.
\newblock Weak type {$(1,1)$} bounds for a class of periodic
  pseudo-differential operators.
\newblock {\em J. Pseudo-Differ. Oper. Appl.}, 5(4):507--515, 2014.

\bibitem{CardonaDelgadoRuzhansky2022Drift}
D.~Cardona, J.~Delgado, and M.~Ruzhansky.
\newblock Drift diffusion equations with fractional diffusion on compact {Lie}
  groups.
\newblock {\em Journal of Evolution Equations}, 22:Article 84, 2022.

\bibitem{CardonaDelgadoRuzhansky2024Hyperbolic}
D.~Cardona, J.~Delgado, and M.~Ruzhansky.
\newblock Well-posedness for a class of pseudo-differential hyperbolic
  equations on the torus.
\newblock {\em Aequationes Mathematicae}, 98:1019--1038, 2024.

\bibitem{CardonaKumar2021Nuclear}
D.~Cardona and V.~Kumar.
\newblock The nuclear trace of periodic vector-valued pseudo-differential
  operators with applications to index theory.
\newblock {\em Mathematische Nachrichten}, 294(9):1657--1683, 2021.

\bibitem{delgado13}
J.~Delgado.
\newblock {$L^p$}-bounds for pseudo-differential operators on the torus.
\newblock In {\em Pseudo-differential operators, generalized functions and
  asymptotics}, volume 231 of {\em Oper. Theory Adv. Appl.}, pages 103--116.
  Birkh\"auser/Springer Basel AG, Basel, 2013.

\bibitem{denk-nau13}
R.~Denk and T.~Nau.
\newblock Discrete {F}ourier multipliers and cylindrical boundary-value
  problems.
\newblock {\em Proc. Roy. Soc. Edinburgh Sect. A}, 143(6):1163--1183, 2013.

\bibitem{FacklerHytonenLindemulder2020Weighted}
S.~Fackler, T.~P. Hyt{\"o}nen, and N.~Lindemulder.
\newblock Weighted estimates for operator-valued {Fourier} multipliers.
\newblock {\em Collectanea Mathematica}, 71(3):511--548, 2020.

\bibitem{favini-giudetti-yakubov11}
A.~Favini, D.~Guidetti, and Y.~Yakubov.
\newblock {Abstract elliptic and parabolic systems with applications to
  problems in cylindrical domains}.
\newblock {\em Advances in Differential Equations}, 16(11/12):1139 -- 1196,
  2011.

\bibitem{Guinazu-Vergara21}
A.~Gui\~naz\'u and V.~Vergara.
\newblock Fundamental solutions and decay rates for evolution problems on the
  torus {$\mathbb {T}^n$}.
\newblock {\em J. Pseudo-Differ. Oper. Appl.}, 12(3):Paper No. 42, 13, 2021.

\bibitem{HNVW16}
T.~Hyt\"onen, J.~v. Neerven, M.~Veraar, and L.~Weis.
\newblock {\em Analysis in {B}anach spaces. {V}olume {I}: {M}artingales and
  {L}ittlewood-{P}aley Theory}, volume~63 of {\em Ergebnisse der Mathematik und
  ihrer Grenzgebiete}.
\newblock Springer, 2016.

\bibitem{mclean91}
W.~McLean.
\newblock Local and global descriptions of periodic pseudodifferential
  operators.
\newblock {\em Math. Nachr.}, 150:151--161, 1991.

\bibitem{melo97}
S.~T. Melo.
\newblock Characterizations of pseudodifferential operators on the circle.
\newblock {\em Proc. Amer. Math. Soc.}, 125(5):1407--1412, 1997.

\bibitem{molahajloo-shahla-wong10}
S.~Molahajloo and M.~W. Wong.
\newblock Ellipticity, {F}redholmness and spectral invariance of
  pseudo-differential operators on {$\mathbb{S}^1$}.
\newblock {\em J. Pseudo-Differ. Oper. Appl.}, 1(2):183--205, 2010.

\bibitem{nau-saal12}
T.~Nau and J.~Saal.
\newblock {$H^\infty$}-calculus for cylindrical boundary value problems.
\newblock {\em Adv. Differential Equations}, 17(7-8):767--800, 2012.

\bibitem{rabinovich11}
V.~S. Rabinovich.
\newblock Pseudodifferential operators with operator valued symbols. {F}redholm
  theory and exponential estimates of solutions.
\newblock {\em Mem. Differ. Equ. Math. Phys.}, 53:127--153, 2011.

\bibitem{ruzhansky-turunen07}
M.~Ruzhansky and V.~Turunen.
\newblock On the {F}ourier analysis of operators on the torus.
\newblock In {\em Modern trends in pseudo-differential operators}, volume 172
  of {\em Oper. Theory Adv. Appl.}, pages 87--105. Birkh\"auser, Basel, 2007.

\bibitem{ruzhansky-turunen09}
M.~Ruzhansky and V.~Turunen.
\newblock On the toroidal quantization of periodic pseudo-differential
  operators.
\newblock {\em Numer. Funct. Anal. Optim.}, 30(9-10):1098--1124, 2009.

\bibitem{ruzhansky-turunen10}
M.~Ruzhansky and V.~Turunen.
\newblock {\em Pseudo-differential operators and symmetries}, volume~2 of {\em
  Pseudo-Differential Operators. Theory and Applications}.
\newblock Birkh\"auser Verlag, Basel, 2010.
\newblock Background analysis and advanced topics.

\bibitem{ruzhansky-turunen10a}
M.~Ruzhansky and V.~Turunen.
\newblock Quantization of pseudo-differential operators on the torus.
\newblock {\em J. Fourier Anal. Appl.}, 16:943--982, 2010.

\bibitem{saranen-wendland87}
J.~Saranen and W.~L. Wendland.
\newblock The {F}ourier series representation of pseudodifferential operators
  on closed curves.
\newblock {\em Complex Variables Theory Appl.}, 8(1-2):55--64, 1987.

\bibitem{Taylor91}
M.~E. Taylor.
\newblock {\em Pseudodifferential operators and nonlinear {PDE}}, volume 100 of
  {\em Progress in Mathematics}.
\newblock Birkh\"auser Boston, Inc., Boston, MA, 1991.

\bibitem{turunen00}
V.~Turunen.
\newblock Commutator characterization of periodic pseudodifferential operators.
\newblock {\em Z. Anal. Anwendungen}, 19(1):95--108, 2000.

\bibitem{turunen-vainikko98}
V.~Turunen and G.~Vainikko.
\newblock On symbol analysis of periodic pseudodifferential operators.
\newblock {\em Z. Anal. Anwendungen}, 17(1):9--22, 1998.

\bibitem{wong11}
M.~W. Wong.
\newblock {\em Discrete {F}ourier analysis}, volume~5 of {\em
  Pseudo-Differential Operators. Theory and Applications}.
\newblock Birkh\"auser/Springer Basel AG, Basel, 2011.

\bibitem{XiaXiong2019Mapping}
R.~Xia and X.~Xiong.
\newblock Mapping properties of operator-valued pseudo-differential operators.
\newblock {\em Journal of Functional Analysis}, 277(9):2918--2980, 2019.

\end{thebibliography}

\end{document}